\documentclass[12pt,a4paper,reqno]{amsart}

\title{Gauge groups and bialgebroids}
\date{November 2021}
\usepackage{amsmath, amssymb, amsthm}
\usepackage[english]{babel}
\usepackage{fullpage}
\usepackage{enumerate}
\usepackage[all]{xy}
\usepackage{graphicx}
\usepackage{mathtools}
\usepackage{mathrsfs}
\usepackage{hyperref}
\usepackage{verbatim}

\author{Xiao Han, Giovanni Landi}
\address[]{\textit{Xiao Han},
SISSA, via Bonomea 265, 34136 Trieste, Italy 
\newline \indent 
Present address: IMPAN, Jana i J\c edrzeja \'Sniadeckich 8, 00-656 Warszawa, Poland}
\email{xhan@impan.pl}
\address[]{\textit{Giovanni Landi},
Universit\`a di Trieste,
Via A. Valerio, 12/1, 34127  Trieste, Italy
\newline \indent
Institute for Geometry and Physics (IGAP) Trieste, Italy 
\newline \indent and INFN, Trieste, Italy}
\email{landi@units.it}
 
\DeclareMathOperator{\Map}{Map}

\DeclareMathOperator{\Aut}{Aut}
  
\DeclareMathOperator{\Hom}{Hom}
 
\DeclareMathOperator{\U}{U} 

\newcommand*{\ii}{\textup i}         
         
\newcommand*{\id}{\textup{id}}     
\newcommand*{\diag}{\textup{diag}}

\newcommand{\SL}{\mathrm{SL}_q(2)}   
\newcommand{\ASL}{{\mathcal O}(\mathrm{SL}_q(2))}   
\newcommand{\Asq}{\O(\mathrm{S}^2_{q})}   
\newcommand{\sq}{\mathrm{S}^2_{q}}

\numberwithin{equation}{section}

\theoremstyle{plain}

\newtheorem{thm}{Theorem}[section]
\newtheorem{lem}[thm]{Lemma}
\newtheorem{prop}[thm]{Proposition}
 
\newtheorem{defi}[thm]{Definition}

\theoremstyle{remark}
\newtheorem{exa}[thm]{Example}
\newtheorem{rem}[thm]{Remark}

\newcommand{\nn}{\nonumber}
\newcommand{\ot}{\otimes}

\newcommand{\beq}{\begin{equation}}
\newcommand{\eeq}{\end{equation}}

\newcommand{\B}{\mathcal{B}}

\newcommand{\cL}{\mathcal{L}}

\newcommand{\C}{\mathcal{C}}
\renewcommand{\O}{\mathcal{O}}

\newcommand{\IR}{\mathbb{R}}
\newcommand{\IC}{\mathbb{C}}
\newcommand{\IZ}{\mathbb{Z}}

\newcommand{\zero}[1]{{#1}{}_{\scriptscriptstyle{(0)}}}
\newcommand{\one}[1]{{#1}{}_{\scriptscriptstyle{(1)}}}
\newcommand{\two}[1]{{#1}{}_{\scriptscriptstyle{(2)}}}
\newcommand{\onet}[1]{{#1}{}_{\scriptscriptstyle{(1')}}}
\newcommand{\twot}[1]{{#1}{}_{\scriptscriptstyle{(2')}}}
\newcommand{\three}[1]{{#1}{}_{\scriptscriptstyle{(3)}}}
\newcommand{\four}[1]{{#1}{}_{\scriptscriptstyle{(4)}}}

\newcommand{\tuno}[1]{{#1}{}{}^{\scriptscriptstyle{<1>}}}
\newcommand{\tdue}[1]{{#1}{}{}^{\scriptscriptstyle{<2>}}}

\begin{document}

\begin{abstract}
We study the Ehresmann--Schauenburg bialgebroid of a noncommutative principal bundle as a quantization of the gauge groupoid of a classical principal bundle. We show that the gauge group of the noncommutative bundle is isomorphic to the group of bisections of the bialgebroid,
and we give a crossed module structure for the bisections and the automorphisms of the bialgebroid. 
Examples include: Galois objects of Taft algebras, a monopole bundle over a quantum sphere
and a not faithfully flat Hopf--Galois extension of commutative algebras. For each of the latter two examples there is in fact a suitable invertible antipode for the bialgebroid making it a Hopf algebroid.
\end{abstract}

\maketitle
\tableofcontents
\parskip = .75 ex

\newpage
\section{Introduction}
The study of groupoids on the one hand and gauge theories on the other hand is important in different areas of mathematics and physics. In particular these subjects meet in the notion of the gauge groupoid of a principal bundle. In view of the considerable amount of recent work on noncommutative principal bundles it is 
desirable to come up with a noncommutative version of groupoids and study their relations to noncommutative principal bundles. For all of this there is a need for a better understanding of bialgebroids. 

In the present paper, having in mind applications to noncommutative gauge theories, we consider the Ehresmann--Schauenburg bialgebroid associated with a noncommutative principal bundle as a quantization of the classical gauge groupoid. Classically, bisections of the gauge groupoid are closely related to gauge transformations.
In parallel with this result  we show that in a rather general context the gauge group of a noncommutative principal bundle is group isomorphic to the group of bisections of the corresponding 
Ehresmann--Schauenburg bialgebroid. To illustrate the theory we work out all the details of the gauge group of the principal bundle and of the bialgebroid with corresponding group of bisections, for the noncommutative $\U(1)$ bundle over the quantum standard Podle\'s 2-sphere, and for a commutative not faithfully flat Hopf--Galois extension obtained in \cite{Be00}  from a particular coaction on the algebra $\O(\mathrm{SL}(2))$. In fact, in each of these two cases there is also an invertible antipode which satisfies the conditions for a Hopf algebroid. In general, for a bialgebroid there is a coproduct and a counit but not an antipode. Here we wish to emphasise one important property that the $\O(\mathrm{SL}(2))$ example shows, that is that at least for a commutative algebra (of coinvariants) the Hopf--Galois extension 
needs not be faithfully flat for our constructions to be well defined and our results to be valid.

Part of the paper deals with Galois objects. A Galois object of a Hopf algebra $H$ is a noncommutative principal bundle over a point in a sense: a Hopf--Galois extension of the ground field $\IC$. In contrast to the classical case where a bundle over a point is trivial, for the isomorphism classes of noncommutative principal bundles over a point this needs not be the case. An antipode can always be defined for the Ehresmann--Schauenburg bialgebroid of a Galois object which (the bialgebroid that is) is then a Hopf algebra. 
Notable examples are group Hopf algebras $\IC[G]$, whose corresponding principal bundle are $\IC[G]$-graded algebras and are classified by the cohomology group $H^2(G, \IC^{\times})$, and Taft algebras $T_{N}$. The equivalence classes of  $T_{N}$-Galois objects are in bijective correspondence with the abelian group $\IC$. 
Thus, part of the paper concerns the Ehresmann--Schauenburg bialgebroid of a Galois object and  corresponding groups of bisections, been they algebra maps from the bialgebroid to the ground field (and thus characters) or more general transformations. For these bialgebroids some of the results we report could be and have been obtained in an abstract and categorical way. Here we re-obtained them in an explicit and more workable fashion, for potential applications to noncommutative  gauge theory.

Automorphisms of a (usual) groupoid  with natural transformations form a strict 2-group or, equivalently, a crossed module.  The crossed module involves the product of bisections and the composition of automorphisms, together with the action of automorphisms on bisections by conjugation. Bisections are the 2-arrows from the identity morphisms to automorphisms, and the composition of bisections can be viewed as the horizontal composition of 2-arrows. In the present paper this construction is extended to the Ehresmann--Schauenburg bialgebroid of a Hopf--Galois extension by constructing a crossed module for the bisections and the automorphisms of the bialgebroid. 

The paper is organised as follows.
After a recap in \S \ref{s2} of algebraic preliminaries and notation, in \S \ref{s3} we give the relevant concepts for noncommutative principal bundles (Hopf--Galois extensions), gauge groups and bialgebroids that we need.
We then work out in \S \ref{ncu1b} the gauge group
for the noncommutative $\U(1)$ principal bundle over 
the quantum sphere and in \S \ref{gtff} for 
a commutative not faithfully flat Hopf--Galois extension associated to $\O(\mathrm{SL}(2))$.
In \S \ref{s4}, we first have Ehresmann--Schauenburg bialgebroids and the group of their bisections. Then we show that the group of gauge transformations of a noncommutative principal bundle  is group isomorphic to the group of bisections of the corresponding Ehresmann--Schauenburg bialgebroid. 
In \S \ref {s5} we describe the Hopf algebroid structure for the $\U(1)$ principal bundle over 
the quantum sphere in \S \ref{ncu1b} and for the commutative 
not faithfully flat Hopf--Galois extension out of $\O(\mathrm{SL}(2))$ considered in \S \ref{gtff}. 
In \S \ref{s6} we consider Galois objects with several examples, such as Galois objects for a cocommutative Hopf algebra, in particular group algebras, regular Galois objects (Hopf algebras as self-Galois objects) and Galois objects of Taft algebras. Finally, in \S \ref{s7}, we study the crossed module (or 2-group) structure coming from the bisections and the automorphism group of a Ehresmann--Schauenburg bialgebroid. When restricting to Hopf algebras one is lead to the representation theory of crossed modules on them. In the present paper we work out this construction for the Taft algebras; more general results will be reported elsewhere.  

\section{Algebraic preliminaries}\label{s2}
We recall here some known facts from algebras and coalgebras and corresponding modules and comodules.
We also recall the more general notions of rings and corings over an algebra as well as the associated notion of bialgebroid.
We move then to Hopf--Galois extensions, as noncommutative principal bundles, and to the definitions of gauge groups.

\subsection{Algebras, coalgebras and all that}\label{acat}
We work over the field $\IC$ of complex numbers but this could be substituted by any commutative field $k$. Algebras (coalgebras) are assumed to be unital and associative (counital and coassociative) with morphisms 
of algebras taken to be unital (of coalgebras taken counital). 
For the coproduct of a coalgebra $\Delta : H \to H \ot H$ we use the Sweedler notation
$\Delta(h)=\one{h}\ot\two{h}$ (sum understood) and
its iterations: $\Delta^n=(\id \ot  \Delta_H) \circ\Delta_H^{n-1}: h \mapsto \one{h}\ot \two{h} \ot
\cdots \ot h_{\scriptscriptstyle{(n+1)\;}}$.
We denote by $*$ the convolution product in the dual vector space
$H':=\mathrm{Hom}(H,\IC)$, $(f * g) (h):=f(\one{h})g(\two{h})$.
The antipode of a Hopf algebra $H$ is denoted $S$.  

Given an algebra $A$, a left $A$-module is a vector space $V$ carrying a left $A$-action, that is with a $\IC$-linear map
$\triangleright_V: A \ot  V \to V$ such that
$$
(a b) \triangleright_V v = a \triangleright_V (b \triangleright_V v )  \, ,\quad 1 \triangleright_V v = v \, .
$$
Dually, with a coalgebra $(H,\Delta)$, a right $H$-comodule is a vector space $V$ carrying a
right $H$-coaction, that is with a $\IC$-linear map $\delta^V : V\to V\ot  H$
such that
$$
(\id\ot  \Delta)\circ \delta^V = (\delta^V\ot  \id)\circ \delta^V  \, ,  \quad
(\id\ot  \varepsilon) \circ \delta^V = \id  \, .
$$
In Sweedler-like notation, $\delta^V (v) = \zero{v}\ot  \one{v}$, and the right $H$-comodule properties
read
\begin{align*}
\zero{v} \ot  \one{(\one{v})}\ot  \two{(\one{v})} = \zero{(\zero{v})} \ot \one{(\zero{v})} \ot  \one{v} =: \zero{v} \ot\one{v} \ot \two{v}  \, , 
\end{align*}
and $\zero{v} \,\varepsilon (\one{v}) = v$, for all $v\in V$.
The $\IC$-vector space tensor product $V\ot  W$ of two $H$-comodules is a $H$-comodule with the right tensor product $H$-coaction
\begin{align}\label{deltaVW}
 \delta^{V\ot  W} :V\ot  W \longrightarrow  V\ot  W\ot  H \, , \quad
 v\ot  w \longmapsto \zero{v}\ot  \zero{w} \ot
 \one{v}\one{w}  \, .
\end{align}
An $H$-comodule map $\psi:V\to W$  between two $H$-comodules is a $\IC$-linear map
$\psi:V\to W$ which is $H$-equivariant (or $H$-colinear), that is $\delta^W\circ \psi=(\psi\ot \id)\circ \delta^V$.

In particular, a right $H$-comodule algebra is an algebra  $A$  which is a right $H$-comodule such that the multiplication and unit of $A$ are morphisms of $H$-comodules.
This is equivalent to requiring the coaction $\delta^A: A\to A\ot  H$ to be
a morphism of unital algebras (where $A\ot  H$ has the usual tensor
product algebra structure).
Corresponding morphisms are $H$-comodule maps which are also algebra maps.

In the same way, a right $H$-comodule coalgebra is a coalgebra $C$ which is a right $H$-comodule and such that the
coproduct and the counit of $C$ are morphisms of $H$-comodules. Explicitly, this means that, for each $c \in C$,
\begin{align*}
\zero{(\one{c})} \ot \zero{(\two{c})} \ot \one{(\one{c})} \one{(\two{c})}
& =
\one{(\zero{c})} \ot \two{(\zero{c})} \ot \one{c} \, ,
\end{align*}
and $\varepsilon(\zero{c}) \one{c} =\varepsilon(c) 1_H$. Corresponding morphisms are $H$-comodule maps which are also coalgebra maps.
There are right $A$-modules and left $H$-comodule versions of the above.

Next, let $H$ be a coalgebra and let $A$ be a right $H$-comodule algebra.
An $(A,H)$-relative Hopf module $V$ is a right $H$-comodule with a compatible left $A$-module structure.
That is  the left action $\triangleright_V: A\ot  V \to V$ is a morphism of $H$-comodules:
$\delta^V \circ \triangleright_V = (\triangleright_V \ot \id) \circ \delta^{A \ot V}$.
Explicitly, for all $a\in A$ and $v\in V$: 
$\zero{(a \triangleright_V v)} \ot \one{(a \triangleright_V v)} = \zero{a} \triangleright_V \zero{v} \ot \one{a}\one{v}.$

A morphism of
$(A,H)$-relative Hopf modules is a morphism of right $H$-comodules  which is also a morphism of left $A$-modules.
In a similar way one can consider the case for the algebra $A$ to be acting on the right, or  with a left and a right 
(bimodule) $A$-actions.

For an algebra $B$ a {\em $B$-ring} is a triple $(A,\mu,\eta)$. Here $A$ is a $B$-bimodule with $B$-bimodule maps
$\mu:A\ot_ {B} A \to A$ and $\eta:B\to A$, satisfying the associativity and unit conditions:
\begin{equation}
\mu\circ(\mu\ot _{B}  \id_A)=\mu\circ (\id_A \ot _{B} \mu), 
\quad 
\mu\circ(\eta \ot _{B} \id_A)=\id_A=\mu\circ (\id_A\ot _{B} \eta).
\end{equation}
A morphism of $B$-rings $f:(A,\mu,\eta)\to (A',\mu',\eta')$ is a
$B$-bimodule map $f:A \to A'$ such that
$f\circ \mu=\mu'\circ(f\ot_{B} f)$ and $f\circ \eta=\eta'$.

From \cite[Lemma 2.2]{Boehm} there is a bijective correspondence between $B$-rings $(A,\mu,\eta)$ and algebra automorphisms
$\eta : B \to A$. Starting with a $B$-ring $(A,\mu,\eta)$, one obtains a multiplication map $A \ot A \to A$ by composing the canonical surjection $A \ot A \to A\ot_B A$ with the map $\mu$. Conversely, starting with an algebra map $\eta : B \to A$, a $B$-bilinear associative multiplication $\mu:A\ot_ {B} A \to A$ is obtained from the universality of the coequaliser $A \ot A \to A\ot_B A$ which identifies an element $ a b \ot a'$ with $ a \ot b a'$.

Dually,  for an algebra $B$ a {\em $B$-coring} is a
triple $(C,\Delta,\varepsilon)$. Here $C$ is a $B$-bimodule with $B$-bimodule maps
$\Delta:C\to C\ot_{B} C$ and $\varepsilon: C\to B$, satisfying the 
coassociativity and counit conditions:
\begin{align}
(\Delta\ot _{B} \id_C)\circ \Delta = (\id_C \ot _{B} \Delta)\circ \Delta, \quad
(\varepsilon \ot _{B} \id_C)\circ \Delta = \id_C =(\id_C \ot _{B} \varepsilon)\circ \Delta.
\end{align}
A morphism of $B$-corings $f:(C,\Delta,\varepsilon)\to
(C',\Delta',\varepsilon')$ is a $B$-bimodule map $f:C \to C'$, such that
$\Delta'\circ f=(f\ot_{B} f)\circ \Delta$ and
$\varepsilon' \circ f =\varepsilon$.

Let $B$ be an algebra.
A {\em left $B$-bialgebroid} $\cL$ consists of a $(B\ot B^{op})$-ring
together with a $B$-coring structures on the same vector space $\cL$ with mutual compatibility conditions
\cite{Take77}.
From what said above, a $(B\ot B^{op})$-ring $\cL$ is the same as an algebra map $\eta : B \ot B^{op} \to \cL$.
Equivalently, one may consider the restrictions
$$
s := \eta ( \, \cdot \, \ot_B 1_B ) : B \to \cL \quad \mbox{and} \quad t := \eta ( 1_B \ot_B  \, \cdot \, ) : B^{op} \to \cL
$$
which are algebra maps with commuting ranges in $\cL$, called the \emph{source} and the \emph{target} map of the
$(B\ot B^{op})$-ring $\cL$. Thus a $(B\ot B^{op})$-ring is the same as a triple $(\cL,s,t)$ with $\cL$ an algebra and $s: B \to \cL$
and $t: B^{op} \to \cL$ both algebra maps with commuting range.

For a left $B$-bialgebroid $\cL$ the compatibility conditions are required to be
\begin{itemize}
\item[(i)] The bimodule structures in the $B$-coring $(\cL,\Delta,\varepsilon)$ are
related to those of the $B\ot B^{op}$-ring $(\cL,s,t)$ via
\begin{equation}\label{eq:rbgd.bimod}
b\triangleright a \triangleleft \tilde{b}:= s(b) t(\tilde{b})a \, \qquad \textrm{for} \,\, b, \tilde{b}\in B, \, a\in \cL.
\end{equation}

\item[(ii)] Considering $\cL$ as a $B$-bimodule as in \eqref{eq:rbgd.bimod},
  the coproduct $\Delta$ corestricts to an algebra map from $\cL$ to
\begin{equation}\label{eq:Tak.prod}
\cL \times_{B} \cL := \left\{\ \sum\nolimits_j a_j\ot_{B} \tilde{a}_j\ |\ \sum\nolimits_j a_jt(b) \ot_{B} \tilde{a}_j =
\sum\nolimits_j a_j \ot_{B}  \tilde{a}_j s(b), \,\,\, \forall  \, b \in B\ \right\},
\end{equation}
where $\cL \times_{B} \cL$ is an algebra via component-wise multiplication.
\\
\item[(iii)] 
The counit $\varepsilon : \cL \to B$ 
satisfies the properties,
\begin{itemize}
\item[(1)] $\varepsilon(1_{\cL})=1_{B}$,   
\item[(2)] $\varepsilon(s(b)a)=b\varepsilon(a) $,   
\item[(3)] $\varepsilon(as(\varepsilon(\tilde{a})))=\varepsilon(a\tilde{a})=\varepsilon(at (\varepsilon(\tilde{a})))$,  \qquad 
for all $b\in B$ and $a,\tilde{a} \in \cL$.
\end{itemize}
\end{itemize}

An automorphism of the left bialgebroid $(\cL, \Delta, \varepsilon,s,t)$ over the algebra $B$ is a pair $(\Phi, \varphi)$
of algebra automorphisms, $\Phi: \cL \to \cL$, $\varphi : B \to B$ such that:
\begin{align}
\Phi\circ s & = s\circ \varphi , \qquad \Phi\circ t = t\circ \varphi ,  \label{amoeba(i)} \\
(\Phi\ot_{B} \Phi)\circ \Delta & = \Delta \circ \Phi , \qquad
\varepsilon\circ \Phi= \varphi\circ \varepsilon \label{amoeba(ii)} .
 \end{align}
In fact, the map $\varphi$ is uniquely determined by $\Phi$ via $\varphi = \varepsilon \circ \Phi \circ s$ 
and one can just say that $\Phi$ is a bialgebroid automorphism. Automorphisms of a bialgebroid $\cL$
form a group $\Aut(\cL)$ by map composition. A \textit{vertical} automorphism is one of the type $(\Phi, \varphi=\id_{B})$.

The pair of algebra maps $(\Phi, \varphi)$ can be viewed as a bialgebroid map (see  \cite[\S 4.1]{kornel}) between two copies of $\cL$ with different source and target maps (and so $B$-bimodule structures). If $s, t$ are the source and target maps on $\cL$, one defines on $\cL$ new source and target maps  
by $s':=s\circ \varphi $ and $t':=t\circ \varphi $ with the new bimodule structure given by $b\triangleright_{\varphi}c\triangleleft_{\varphi}\tilde{b}:=s'(b)t'(\tilde{b})a$, for any $b, \tilde{b}\in B$ and $a\in \cL$ (see \eqref{eq:rbgd.bimod}).
Therefore one gets a new left bialgebroid with product, unit, coproduct and counit not changed.

From the conditions \eqref{amoeba(i)}, $\Phi$ is a $B$-bimodule map: 
$\Phi(b \triangleright c\triangleleft \tilde{b})=b \triangleright_{\varphi}\Phi(c)\triangleleft_\varphi \tilde{b}$. 
The first condition \eqref{amoeba(ii)} is well defined once the conditions \eqref{amoeba(i)} are satisfied (the balanced tensor product is induced by $s'$ and $t'$). Conditions \eqref{amoeba(i)} imply that $\Phi$ is a coring map; 
therefore $(\Phi, \varphi)$ is an isomorphism between the starting bialgebroid and the new one.
  
Finally, we recall from \cite[Def. 4.1]{BS04} the conditions for a Hopf algebroid with invertible antipode. 
 Given a left bialgebroid $(\cL, \Delta, \varepsilon,s,t)$ over the algebra $B$, an invertible antipode
$S : \cL  \to  \cL $ in an algebra anti-homomorphism with inverse $S^{-1} : \cL  \to  \cL $ such that
\beq\label{hopbroid1}
S \circ t = s
\eeq
and satisfying compatibility conditions with the coproduct:
\begin{align}\label{hopbroid2}
\onet {(S^{-1}\two{h})} \ot_B \twot {(S^{-1}\two{h})}\one{h} & = S^{-1} h \ot_B 1_\cL \nn \\
\onet {(S\one{h})} \two{h} \ot_B \twot {S(\one{h})} & = 1_\cL \ot_B S h ,
\end{align}
for any $h\in \cL $. These then imply  
$
S(\one{h}) \, \two{h} = t \circ \varepsilon \circ S h. 
$
 
\section{Noncommutative principal bundles}\label{s3}
We start with a brief recall of Hopf--Galois extensions
as noncommutative principal bundles. Then we consider gauge transformations 
as equivariant automorphisms of the total space algebra which are vertical so that they leave invariant the base space algebra.

\subsection{Hopf--Galois extensions}\label{sec:hge}
These extensions are $H$-comodule algebras
$A$ with a canonically defined map $\chi: A\ot_B A\to A\ot H$
which is required to be invertible \cite{HJSc90}.

\begin{defi} \label{def:hg}
Let $H$ be a Hopf algebra and let $A$ be a $H$-comodule algebra with coaction $\delta^A$.
Consider the subalgebra $B:= A^{coH}=\big\{b\in A ~|~ \delta^A (b) = b \ot 1_H \big\} \subseteq 
A$ of coinvariant elements 
with balanced tensor product $A \ot_B A$.
The extension $B\subseteq A$ is called a $H$-\textup{Hopf--Galois extension} if the
\textit{canonical Galois map}
\begin{align*}  
\chi := (m \ot \id) \circ (\id \ot _B \delta^A ) : 
A \ot _B A \longrightarrow A \ot H ,
\quad a' \ot_B a  &\mapsto a' \zero{a} \ot \one{a}
\end{align*}
is an isomorphism.
\end{defi}
\begin{rem}\label{fun-rem}
For a Hopf--Galois extension $B\subseteq A$, we take (apart from Sects.~\ref{gtff} and \ref{biff})
the algebra $A$ to be \emph{faithfully flat} as a right $B$-module.
One possible way to state this property is that for any left $B$-module map $F : M \to N$, 
the map $F$ is injective if and only if the map $\id_A \ot_B F : A  \ot_B M \to A \ot_B N$ 
is injective; injectivity of $F$ implying the injectivity of $\id_A \ot_B F$ would state that $A$ is flat as a right $B-$module (see \cite[Chap. 13]{Wa79}.
%
For a faithfully flat $H$-Hopf--Galois extension 
the category of $(A,H)$-relative Hopf modules is equivalent to the category of left $B$-modules by 
$M \to A\ot_B M$ \cite[Thm.~1]{HJSc90}. 
\end{rem}

The canonical map $\chi$ is a morphism of relative Hopf modules for $A$-bimodules and right $H$-comodules \cite[\S 1.1]{HJSc90b}. Both $A \ot _B A$ and $A \ot H$ are $A$-bimodules .
The left $A$-module structures are left multiplication on the first factor while the right $A$-actions are
$$
(a \ot_B a')a'':=  a \ot_B a'a'' \quad \mbox{and} \quad (a \ot h) a' := a \zero{a'} \ot h \one{a'} \, .
$$
For the $H$-comodule structure, the right tensor product $H$-coaction as in \eqref{deltaVW}:
\beq\label{AAcoact}
\delta^{A\ot  A}: A\ot  A\to A\ot  A\ot  H, \quad a\ot  a'  \mapsto  \zero{a}\ot  \zero{a'} \ot   \one{a}\one{a'} \, ,
\eeq
for all $a,a'\in A$, descends to the quotient $A\ot_B A$ because $B\subseteq A$ is the
subalgebra of $H$-coinvariants.
Similarly,  $A \ot H$ is endowed with the tensor product coaction, where   one regards the Hopf algebra $H$ as a right $H$-comodule
with the right adjoint $H$-coaction
$$
\mathrm{Ad} :
h \longmapsto \two{h}\ot  S(\one{h})\,\three{h}  \, .
$$
The right $H$-coaction on $A \ot H$
is then given, for all $ a\in A, \, h \in H$, by
$$
\delta^{A\ot  H}(a\ot  h) = \zero{a}\ot  \two{h} \ot  \one{a}\,S(\one{h})\, \three{h} \in A \ot H \ot H  \, .
$$

Since the canonical Galois map $\chi$ is left $A$-linear, its inverse is
determined by the restriction $\tau:=\chi^{-1}_{|_{1_A \ot H}}$, named \textit{translation map},
\begin{eqnarray*}
\tau =\chi^{-1}_{|_{1_A \ot H}} :  H\to A\ot _B A ~ ,
\quad h \mapsto \tau(h) = \tuno{h} \ot_B \tdue{h} \, .
\end{eqnarray*}
Thus by definition:
\beq
\label{p7}
\tuno{h}\zero{\tdue{h}}\ot \one{\tdue{h}} = 1_{A} \ot h  \, .
\eeq
The translation map enjoys a number of properties \cite[3.4]{HJSc90b}
that we list here for later use. 
For any $h, k \in H$ and $a\in A$, $b\in B$:  
\begin{align}
\tuno{h} \ot_B \zero{\tdue{h}} \ot \one{\tdue{h}} &= \tuno{\one{h}} \ot_B \tdue{\one{h}} \ot \two{h} \label{p4}  \, , \\
\zero{\tuno{h}} \ot_B {\tdue{h}} \ot \one{\tuno{h}} &= \tuno{\two{h}}  \ot_B \tdue{\two{h}} \ot S(\one{h})  \, , \label{p1}
\end{align}
\begin{align}
\label{p5}
\tuno{h}\tdue{h} &= \varepsilon(h)1_{A}  \, , \\
\label{p3}
\zero{a}\tuno{\one{a}}\ot_{B}\tdue{\one{a}} &= 1_{A} \ot_{B}a  , 
\end{align}
\begin{align}
\label{p2}
\tuno{(h k)}\ot_{B}\tdue{(h k)} &= \tuno{k}\tuno{h}\ot_{B}\tdue{h}\tdue{k}  \, , \\ 
\label{p6}
\tuno{\one{h}}\ot_{B}\tdue{\one{h}}\tuno{\two{h}}\ot_{B}\tdue{\two{h}} & 
=\tuno{h}\ot_{B}1_{A}\ot_{B}\tdue{h}  \, , 
\end{align}
\begin{align}
\label{p8}
b\, \tuno{h} \ot_B \tdue{h} &= \tuno{h} \ot_B \tdue{h} \, b \,. 
\end{align}

Two Hopf--Galois extensions $A,A' $ for the Hopf algebra $H$ of the
algebra $B$ are isomorphic provided there exists an isomorphism of
$H$-comodule algebras $A \to A'$. This is the algebraic counterpart
for noncommutative principal bundles of the geometric
notion of isomorphism of principal $G$-bundles over the same base space.

\subsection{The group of gauge transformations}
In \cite{brz-tr} gauge transformations for a noncommutative principal bundles were defined to be
invertible and unital comodule maps, with no additional requirement. In particular they were not
asked to be algebra maps. However, the resulting gauge group might be
very big, even in the classical case. For example the gauge group of a $G$-bundle over a point would be
much bigger than the structure group $G$. In contrast, in \cite{PGC20} gauge transformations were
taken to be algebra homomorphisms. This property implies in particular that they are invertible.

\begin{prop}\label{prop:gsogg}
Let $B=A^{coH}\subseteq A$ be a faithfully flat Hopf--Galois  extension. Then the collection $\Aut_H(A)$
of right $H$-comodule unital algebra maps of $A$ into itself which restrict to the identity on the subalgebra $B$
is a group for map composition. 

Moreover, recall the notation $\tau(h) = \tuno{h} \ot_B \tdue{h}$ for the translation map. 
Then for $F\in\Aut_H(A)$ its inverse $F^{-1} \in \Aut_H(A)$ is given, for all $a\in A$, by
\beq\label{inv-F2} 
F^{-1}(a) = \zero{a} F (\tuno{\one{a}}) \, \tdue{\one{a}} \,  .
\eeq
\end{prop}
\begin{proof}
That vertical $H$-comodule algebra maps are invertible is in \cite[Rem. 3.11]{HJSc90}. 
We check the expression of the inverse in \eqref{inv-F2}. 
The $B$-linearity and the algebra map property assure that the inverse is well defined: $
\zero{a} F (\tuno{\one{a}} b) \, \tdue{\one{a}} = \zero{a} F (\tuno{\one{a}}) \, b \tdue{\one{a}}$, 
for $b \in B$; also $F^{-1}(b) = b$ for $b \in B$. 
For any $a\in A$, using the $H$-equivariance of $F$, 
\begin{align}
F^{-1}(F(a))& = \zero{F(a)}F(\tuno{\one{F(a)}})  \, \tdue{\one{F(a)}}= F(\zero{a})F(\tuno{\one{a}}) \, \tdue{\one{a}} \nn \\ 
& = F(\zero{a}\tuno{\one{a}}) \, \tdue{\one{a}} = a. \label{f1f}
\end{align}
Now, the action of the canonical map $\chi$ yields
the equality
\beq\label{chi-eq}
1_A  \ot_B \zero{a} F (\tuno{\one{a}}) \tdue{\one{a}} = \zero{a} F (\tuno{\one{a}}) \ot_B \tdue{\one{a}}.
\eeq
Indeed, using \eqref{p4}, for the right hands side:
\begin{align*}
\chi (\zero{a} F (\tuno{\one{a}}) \ot_B \tdue{\one{a}}) & = 
\zero{a} F (\tuno{\one{a}}) \, \zero{\tdue{\one{a}}} \ot \one{\tdue{\one{a}}} \\
& = \zero{a} F (\tuno{\one{\one{a}}}) \, \tdue{\one{\one{a}}} \ot \two{\one{a}} \\ 
& = \zero{a} F (\tuno{\one{a}}) \, \tdue{\one{a}} \ot \two{a} .
\end{align*}
Next, using equivariance, \eqref{p4} and \eqref{p1}, for the left hand side:
\begin{align*}
\chi (1_A  \ot_B \zero{a} F (\tuno{\one{a}}) \tdue{\one{a}} )
& = \zero{a} \zero{F (\tuno{\two{a}})} \, \zero{\tdue{\two{a}}} \, \ot 
\one{a} \one{F (\tuno{\two{a}})} \, \one{\tdue{\two{a}}} \\ 
& = \zero{a} F (\zero{\tuno{\two{a}}}) \, \zero{\tdue{\two{a}}} \, \ot 
\one{a} \one{\tuno{\two{a}}} \, \one{\tdue{\two{a}}} \\ 
& = \zero{a} F (\tuno{\three{a}}) \, \zero{\tdue{\three{a}}} \, \ot 
\one{a} S(\two{a}) \, \one{\tdue{\three{a}}} \\ 
& = \zero{a} F (\tuno{\one{a}}) \, \zero{\tdue{\one{a}}} \ot \one{\tdue{\one{a}}} \\
& = \zero{a} F (\tuno{\one{a}}) \, \tdue{\one{a}} \ot \two{a} .
\end{align*}
Then bijectivity of the canonical map $\chi$ yields the identity \eqref{chi-eq}. 
With $B$-linearity of $F$, using $\id_A \ot_B F$ on both sides
of equality \eqref{chi-eq}, 
the right faithful flatness
leads to
\begin{align*}
F(F^{-1}(a))& = F \big(\zero{a}F(\tuno{\one{a}}) \big) F(\tdue{\one{a}}) = \zero{a}F(\tuno{\one{a}})F(\tdue{\one{a}})\\
& = \zero{a}F(\tuno{\one{a}} \, \tdue{\one{a}}) = a.
\end{align*}
Thus $F^{-1}$ is the inverse map of $F \in \Aut_H(A)$. That $F^{-1}$ is a $H$-comodule algebra map follows directly from such properties of $F$. 
\end{proof}

\noindent
Elements $F\in\Aut_H(A)$ preserve the (co)-action of the structure quantum group since they are such that
$\delta^A \circ F = (F \ot \id) \delta^A$ (or $\zero{F(a)} \ot \one{F(a)} = F(\zero{a}) \ot \one{a}$). 
And they also preserve the base space algebra $B$.
This group will be called the \emph{gauge group}.

\begin{rem}
A similar proposition was given in \cite{PGC20}, for $H$ a coquasitriangular  Hopf  algebra, and $A$ a quasi-commutative $H$-comodule algebra. As a consequence, $B$ is in the centre of $A$. In the present paper, there is no restriction on the coinvariant subalgebra B.
 \end{rem}

\subsection{Noncommutative $\U(1)$-bundles}\label{ncu1b}
Let $G$ be a  group and $\IC[G]$ be its group algebra.
Its elements are finite sums $\sum \lambda_g \, g$ with $\lambda_g \in \IC$.
The algebra product follows from the group product in $G$, with unit $1_{\IC[G]}=e$,
the neutral element of $G$. The coproduct, counit and antipode, making $\IC[G]$ a Hopf algebra, are 
$\Delta(g) = g \ot g$, $\varepsilon(g)=1$, $S(g)=g^{-1}$.

It is known that $\IC[G]$-Hopf--Galois extensions are the same as strongly graded algebras over $G$.  Now, an algebra $A$ is $G$-graded, $A=\oplus_{g \in G} \, A_g$ with $A_g A_{h} \subseteq A_{gh}$ for all $g,h \in G$,
if and only if $A$ is a right $\IC[G]$-comodule  algebra with coaction
$\delta^A: A \to A \ot \IC[G]$ given by $a \mapsto \sum a_g \ot g$ for $a= \sum a_g$, $a_g \in A_g$.
Then the algebra $A$ is \emph{strongly} $G$-graded, that is $A_g A_h =A_{gh}$,
if and only if $B=A_e:=A^{co\, \IC[G]} \subseteq A$ is Hopf--Galois (see \cite[Thm.8.1.7]{mont}).

Let us concentrate on taking 
$H=\O(\U(1)) := \IC[z,z^{-1}] / (1 - z z^{-1})$
where $(1 - z z^{-1})$ is the ideal generated by $1 - z z^{-1}$ in the polynomial algebra $\IC[z,z^{-1}]$ in two variables. The Hopf algebra structure of $H$ is now, for all $n \in \IZ$,  
the coproduct $\Delta : z^n \mapsto z^n \ot z^n$, the antipode $S : z^n \mapsto z^{-n}$ 
and the counit $\varepsilon : z^n \mapsto 1$. 

A \emph{strongly} graded $\IZ$-algebra $A = \oplus_{n \in \IZ} \, A_{n}$ with coaction determined by 
$$
\delta^A(a): A \to A \ot \O(\U(1)) , \qquad x \mapsto x \ot z^{-n} \, , \, \, \, \textup{for} \, \, \, x \in A_{n} \, 
$$
results into a right comodule algebra for $H=\O(\U(1))$ which will be referred to as a noncommutative $\U(1)$ principal bundle over the algebra $B:=A_0$. 

From \cite[Cor.~I.3.3]{NaVO82}, for $A$ a strongly $\IZ$-graded algebra,
the right-modules $A_{1}$ and $A_{-1}$ are finitely generated and projective over $A_{0}$.
In fact, the total space algebra can be recovered out of the `line bundles' $A_{1}$ and $A_{-1}$ as a Pimsner algebra \cite{AKL}.

Then, if $F\in\Aut_H(A)$ is a gauge transformation the equivariance 
$\zero{F(a)} \ot \one{F(a)} = F(\zero{a}) \ot \one{a}$ implies that  $F$ respects the grading 
and in fact $F$ is completely determined by its restrictions to $A_{1}$ and $A_{-1}$ 
as $B$-module maps, given that $F$ is required to be the identity on $B$ and that it is then extended as an algebra map. 

Let us consider an explicit example of the above construction, that is the $\U(1)$ principal bundle over the standard Podle\'s sphere $\sq$ of \cite{Po87}.  
With $q\in\IR$ a deformation parameter,  
the coordinate algebra $\ASL$ of the quantum group $\SL$ is the 
algebra generated by elements $a$, $c$ and $d $, $b$ with relations
\begin{align}\label{reltot}
a\,c &= q\,c\,a \quad \mbox{and} \quad b \,d  =q\,d \, b  ,      
\quad a\, b =q\, b\,a \ 
\quad \mbox{and} \quad c\,d  =q\,d \,c , \nn  \\
c\, b &= b \,c \, , \quad a \, d  - d \, a  = (q - q^{-1}) \,b \, c  
\quad \mbox{and} \quad  d \, a - q^{-1} b\, c = 1 . 
\end{align}
Then the Hopf algebra $H=\O(\U(1))$ coacts on the algebra $\ASL$ via   
\begin{align}\label{rco}
\delta(a) = a \ot z   ,  \quad  
\delta(d ) = d \ot z^{-1} \, \quad \mbox{and}& \quad
 \delta(c) = c \ot z   , \quad  
\delta(b) = b \ot z^{-1} .  
\end{align}
The subalgebra of coinvariant elements in $\ASL$ for this coaction is the coordinate algebra $B=\O(\sq)$ of the standard Podle\'s sphere
 $ 
\Asq := \ASL^{\U(1)}
$. 
As a set of generators for $\O(\sq)$ one may take 
\beq\label{genB}
B_{-} := - q^{-1} a\, b ,  \qquad 
B_{+} := c\,d  \quad \mbox{and} \quad 
B_{0} := - q^{-1} c\, b, 
\eeq
for which one finds the relations
\begin{align} \label{relbas}
B_{-}\,B_{0} &= q^{2}\, B_{0}\,B_{-} \quad \mbox{and} \quad
B_{+}\,B_{0} = q^{-2}\, B_{0}\,B_{+} ,  \nn \\
B_{-}\,B_{+} &= q^2\, B_{0} \,\big( 1 - q^2\, B_{0} \big) \quad
\mbox{and} \quad B_{+}\,B_{-}= B_{0} \,\big( 1 - B_{0} \big) .
\end{align}
The algebra inclusion $\Asq \subset \ASL$, a noncommutative principal
bundle \cite{BM93}, is a faithfully flat Hopf--Galois extension. The translation map on generators of $\O(\U(1))$ is
\beq\label{trau1}
\tau(z) = d \ot_B a - q^{-1} b \ot_B c , \qquad \tau(z^{-1}) = a \ot_B d - q c \ot_B b .
\eeq
The total space algebra decomposes as $\ASL=\bigoplus\nolimits_{n\in\IZ}\, A_n$  
where
\beq\label{libu} 
A_n := \big\{x \in \ASL ~ \big|~\delta(x) = x \ot z^{-n} \big\} .
\eeq 
In particular as $B$-modules, $A_{-1}$ is generated  by $a,c$ while $A_{1}$ is generated 
by $d, b$. Any gauge transformation will be then determined by the images 
\begin{align}\label{fun}
F(a) & = X a + Y c , \qquad F(c) = Z a + W c , \nn \\
F(d ) & = \widetilde X  d + \widetilde Y b , \qquad F(b) = \widetilde Z d + \widetilde W b ,
\end{align}
with coefficients which are elements in the algebra $B$, and extended as an algebra map.

Let us first consider the classical case, $q=1$ of commutative algebras, to clarify the structures.  Asking for the coinvariant generators in \eqref{genB} to be left unchanged by $F$ in \eqref{fun} reduces the coefficient to a single one $X$, any non vanishing function from $\mathrm{S}^2 \to \IC$:
\beq\label{qgt}
F(a) = X a , \qquad F(c) = X c ,  \qquad \mbox{and} \qquad F(d ) =  X^{-1} d , \qquad F(b) = X^{-1} b , 
\eeq
and the sphere relation is automatically satisfied. We get 
$\Aut_H \mathrm{SL}(2) = \Map(\mathrm{S}^2 \to \IC^*)$. 

In contrast, when $q\not=1$,  requiring that $F$ be an algebra map and so to respect the commutation relations in \eqref{reltot}, one gets 
$X\in \IC^*$ since the centre of $\Asq$ is just the algebra $\IC$. Thus $\Aut_H \mathrm{SL}_q(2) = \IC^*$, the non vanishing complex numbers. 


\subsection{A gauge group without faithfully flatness}\label{gtff}
We give an example of the above construction of the gauge group for a Hopf--Galois extension 
over a commutative algebra of coinvariants 
which is flat but not faithfully flat. 

This example was studied in \cite[Ex.~2.4]{Be00}.
Consider the Hopf algebra $H=\IC[x]$ with $x$ a primitive element. 
Let $A = \O(\mathrm{SL}(2)) = \IC[a, b, c, d]/(ad - bc -1)$. 
The algebra $A$ is made into a $H$-comodule algebra with coaction $\delta : A \to A \ot H$ given on generators by
\beq\label{conff}
\delta(a) = a \ot 1 + c \ot x, \quad \delta(b) = b \ot 1 + d \ot x, \quad \delta(c) = c \ot 1, \quad
\delta(d) = d \ot 1 .
\eeq
Then, the algebra of coinvariants is $B= \IC[c, d]$ and the inclusion 
$B = A^{coH} \subset A$ is a Hopf--Galois extension (that is the corresponding canonical map is bijective).
It is shown in \cite[Ex.~2.4]{Be00},  that the extension is flat but not faithfully flat. 
It is easy to see that the corresponding translation map $\tau : H \to A \ot_B A$ is given by
\beq\label{invnff}
\tau(1) = 1\ot_B 1 , \quad \tau(x) = \tuno{x} \ot_B \tdue{x} = a \ot_B b - b \ot_B a  .
\eeq

 Consider then the group $\Aut_H(A)$ of gauge transformations. Any such a map 
is determined by its values on the generators, being $F(c) = c$ and $F(d) = d$ by the $B$-linearity. We claim that given $F(a)$ and $F(b)$ their inverse is as in the formula \eqref{inv-F2}: 
\begin{align*}
F^{-1}(a) = a + c \, \big(F(a) b - F(b) a \big) , \quad F^{-1}(b) = b + d \, \big(F(a) b - F(b) a \big) ,
\end{align*}
clearly together with $F^{-1}(c)=c$ and $F^{-1}(d)=d$.
\begin{lem}\label{lecb}
Let $F$ be any gauge transformation $: A \to A$. Then,
$$
F(a) b - F(b) a \in B.
$$
\end{lem}
\begin{proof}
A direct computation using the equivariance:
\begin{align*}
\delta \big(F(a) b & - F(b) a \big) = 
(F(\zero{a}) \ot \one{a}) \, \delta(b) - (F(\zero{b}) \ot \one{b}) \, \delta(a) \\
& = \big( F(a) b - F(b) a \big) \ot 1 + \big(F(a) d - F(b) c + c b -d a\big) \ot x + (cd - dc) \ot x^2 \\
& = \big( F(a) b - F(b) a \big) \ot 1 , 
\end{align*}
and also using the `commutation relations' of the algebra (just commutativity in the present case) to remove the terms in $x$ and $x^2$.
\end{proof}
\noindent
Next one checks that \eqref{p3} is satisfied for all generators. 
Then, $F^{-1} \circ F = \id_A$ goes as in \eqref{f1f} given the expressions in \eqref{invnff}:
using equivariance and $B$-linearity,
\begin{align*}
F^{-1}(F(a))& = F(a) + c \, \big(F(a) b - F(b) a \big)  \\
& = F(a) + F(a) (a d - 1) - F( ad - 1) a \\
& = a + F(a) a d - F(a) d a = a , 
\end{align*}
with the relation $ad - bc = 1$. Conversely, from Lemma \ref{lecb}, using $B$-linearity and the determinant condition: 
\begin{align*}
F(F^{-1}(a)) &= F(a) + c \big (F(a) b - F(b) a \big) = F(a) + c F(a) b - F(cb) a \\
& = F(a) + c F(a) b + a - F(ad) a = a + F(a) + F(a) (cb - d a ) \\
& = a + F(a) - F(a) = a .
\end{align*}
A similar computation shows that $F^{-1}(F(b)) = b = F(F^{-1}(b))$. The group $\Aut_H(A)$ is not trivial. 
Besides the identity map it contains for instance unital maps of the kind  
\beq \label{gu1b}
F(a) = a + h \, c , \qquad F(b) = b + h \, d 
\eeq
for $h$ an arbitrary element in $B$, and $F(c)=c$, $F(d)=d$, extended as an algebra map. 
This $F$ is equivariant and preserves the determinant  condition: $F(ad-bc) = F(1)=1$. 

\section{Ehresmann--Schauenburg bialgebroids}\label{s4}
To any Hopf--Galois extension $B=A^{co \, H}\subseteq A$ one associates a $B$-coring  and a bialgebroid \cite{schau2} (see \cite[\S 34.13 and 34.14]{BW}).
These can be viewed as a quantization of the gauge or Ehresmann groupoid that is associated to a principal fibre bundle (see  \cite{KirillMackenzie}).

\subsection{The Ehresmann coring} 
The coring can be given in a few equivalent ways.
Let $B=A^{co \, H}\subseteq A$ be a 
Hopf--Galois extension with right coaction $\delta^{A} : A \to A \ot H$.
Recall the diagonal coaction \eqref{AAcoact}, given for all $a, a' \in A$, by
$$
\delta^{A\ot  A}: A\ot  A\to A\ot  A\ot  H, \quad a\ot  a'  \mapsto
\zero{a}\ot  \zero{a'} \ot   \one{a}\one{a'} \, .
$$
Let $\tau$ be the translation map of the Hopf--Galois extension. We have the following: 		
\begin{lem}
The $B$-bimodule of coinvariant elements for the diagonal coaction, 
\beq  \label{ec2}
(A\ot A)^{coH} = \{a\ot  \tilde{a}\in A\ot  A \, ; \,\, \zero{a}\ot  \zero{\tilde{a}}\ot  \one{a}\one{\tilde{a}}=a\ot  \tilde{a}\ot  1_H \}
\eeq
is the same as the $B$-bimodule
\begin{equation}\label{ec1}
\C  :=\{a\ot  \tilde{a}\in A\ot  A \, ; \,\, \zero{a}\ot  \tau(\one{a})\tilde{a}=a\ot  \tilde{a}\ot _B 1_A\}.
\end{equation}
\end{lem}
\begin{proof}
The $B$-bimodule structure of $\C$ is left and right multiplication by elements of $B$.
Let $a\ot \tilde{a}\in (A\ot A)^{coH}$. By applying $(\id_{A}\ot \chi)$ on $\zero{a}\ot\tuno{\one{a}}\ot_{B}\tdue{\one{a}}\tilde{a}$, we get
\begin{align*}
    \zero{a}\ot \tuno{\one{a}}\zero{\tdue{\one{a}}}\zero{\tilde{a}}\ot \one{\tdue{\one{a}}}\one{\tilde{a}} & =\zero{a}\ot \zero{\tilde{a}}\ot \one{a}\one{\tilde{a}} \\ &=a\ot \tilde{a}\ot 1_H =a\ot \chi(\tilde{a}\ot_{B} 1_A) \\
    & = (\id_{A}\ot \chi) (a \ot \tilde{a} \ot_{B} 1_A),
\end{align*}
where the first step uses \eqref{p7}. This shows that $(A\ot A)^{coH}\subseteq \C$.

\noindent
Conversely, let $a\ot \tilde{a}\in \C$. By applying $(\id_{A}\ot \chi^{-1})$ on $\zero{a}\ot\zero{\tilde{a}}\ot \one{a}\one{\tilde{a}}$ and using the fact that $\chi^{-1}$ is left $A$-linear and \eqref{p2}, we get
\begin{align*}
    \zero{a}\ot \zero{\tilde{a}} \tuno{\one{\tilde{a}}}\tuno{\one{a}}\ot_{B}\tdue{\one{a}}\tdue{\one{\tilde{a}}}&=\zero{a}\ot\tuno{\one{a}}\ot_{B}\tdue{\one{a}}\tilde{a}=a\ot \tilde{a}\ot_{B} 1_A \\ &=(\id_{A}\ot \chi^{-1})(a\ot \tilde{a}\ot 1_H),
\end{align*}
where in the first step we used  \eqref{p3}. This shows that $ \C \subseteq (A\ot A)^{coH}$. \end{proof}

We have then the following definition \cite{schau2} (see \cite[\S 34.13]{BW}).
\begin{defi}\label{def:ec}
Let $B=A^{co \, H}\subseteq A$ be a 
faithfully flat 
Hopf--Galois extension with translation map $\tau$.
Then the $B$-bimodule $\C $ in \eqref{ec1} is a $B$-coring with %
coproduct and counit:
\begin{align}\label{copro}
\Delta(a\ot  \tilde{a}) = \zero{a}\ot  \tau(\one{a})\ot \tilde{a} 
& = \zero{a} \ot \tuno{\one{a}} \ot_B \tdue{\one{a}} \ot \tilde{a} , \\
\varepsilon(a\ot  \tilde{a}) &= a\tilde{a} . \label{counit}
\end{align}
\end{defi}
\noindent
Applying the map $m_A\ot \id_H$ to elements of \eqref{ec2} one gets  $a\tilde{a}\in B$.
The above $B$-coring is called the \textit{Ehresmann} or \textit{gauge coring}; we denote it $\C(A, H)$.
Using the known relation between the coinvariants of a tensor product of comodules and their 
cotensor product \cite[Lemma 3.1]{HJSc90b}, the coring $\C(A, H)$ can also be given as a 
cotensor product $A\, \square \, {}^H\!\!A$.

The Ehresmann coring of a Hopf--Galois extension is
in fact a bialgebroid \cite{schau2}, called the \textit{Ehresmann--Schauenburg bialgebroid} (see  \cite[34.14]{BW}).
One see that $\C(A, H) = (A\ot A)^{coH}$ is a subalgebra of $A \ot A^{op}$; indeed, given
$x \ot  \tilde{x}, \, y \ot  \tilde{y} \in (A\ot A)^{coH}$, one computes
$\delta^{A\ot  A}(x y \ot \tilde{y} \tilde{x}) =
\zero{x} \zero{y} \ot \zero{\tilde{y}} \zero{\tilde{x}}  \ot \one{x}\one{y} \one{\tilde{y}} \one{\tilde{x}} =
\zero{x} y \ot \tilde{y} \zero{\tilde{x}}  \ot \one{x} \one{\tilde{x}} =
x y \ot \tilde{y} \tilde{x} \ot  1_H$.

\begin{defi}\label{def:reb}
Let $\C(A, H)$ be the coring associated with a faithfully flat Hopf--Galois extension
$B=A^{co \, H}\subseteq A$. Then $\C(A, H)$
is a (left) $B$-bialgebroid with product
$$
(x \ot  \tilde{x}) \bullet_{\C(A, H)} ({y}\ot  \tilde{y}) = x y \ot \tilde{y} \tilde{x} ,
$$
for all $x \ot  \tilde{x}, \, y \ot  \tilde{y} \in \C(A, H)$ (and unit $1_A\ot  1_A$). The target and the  source  maps are
$$
t(b)=1_A \ot  b \quad \textup{and} \quad s(b)=b\ot 1_A.
$$
\end{defi}
\noindent
We refer to \cite[34.14]{BW} for the check that all defining properties are satisfied. When there is no risk of confusion we drop the decoration $\bullet_{\C(A, H)}$ in the product. 

\subsection{Bisections and gauge groups}
The bialgebroid of a Hopf--Galois extension can be viewed as a quantization (of the dualization)
of the classical gauge groupoid, recalled in Appendix \ref{gaugegroupoid}, of a (classical) principal bundle.
Dually to the notion of a bisection on the classical gauge groupoid there is the notion of a bisection on the
Ehresmann--Schauenburg bialgebroid. 
These bisections correspond to gauge transformations.

The notion of a bisection as in the following definition could be given for any bialgebroid, not only for the Ehresmann--Schauenburg bialgebroid. However, 
for the general case one would need some additional requirements
so to get a proper composition of bisections extending \eqref{mulbis1} below. We shall address this general definition elsewhere. 

\begin{defi}\label{def:bisection} 
Let $\C(A, H)$ be the Ehresmann--Schauenburg bialgebroid of a faithfully flat  
Hopf--Galois extension
$B=A^{coH}\subseteq A$. A bisection of $\C(A, H)$ is a $B$-bilinear unital left character
on the $B$-ring $(\C(A, H), s)$. That is, a map 
$\sigma: \C(A, H)\to B$ such that:
\begin{itemize}
    \item [(1)] $\sigma(1_A \ot 1_A)=1_{B}$, \qquad \textup{unitality, } \item []
    \item[(2)] $\sigma\big(s(b)t(\tilde{b}) (x \ot \tilde{x}) \big)=b\sigma(x \ot \tilde{x}) \tilde{b} $,  
    \qquad $B$-\textup{bilinearity,}  \item []
    \item[(3)] $\sigma\big( (x \ot \tilde{x}) \, s(\sigma(y \ot \tilde{y}))\big)=\sigma\big((x \ot \tilde{x})
    (y \ot \tilde{y}) \big)$, \qquad \textup{associativity, }
\end{itemize}
for all $b, \tilde{b} \in B$ and $x \ot \tilde{x}, \, y \ot \tilde{y}\in \C(A, H)$.
\end{defi}
\noindent
There is also a middle $B$-linearity, that is for any bisection $\sigma$:
\begin{align}\label{equ.bis}
    \sigma(x b \otimes \tilde{x})=\sigma(x \otimes b \tilde{x}),
\end{align}
for any $x \otimes \tilde{x} \in \C(A, H)$ and $b\in B$. Indeed,
\begin{align*}
    \sigma(x\otimes b\tilde{x}) & =\sigma((x \otimes \tilde{x})(1\otimes b))=\sigma((x\otimes \tilde{x})s(\sigma(1\otimes b)))=\sigma((x \otimes \tilde{x})s(b)) =\sigma(x b\otimes \tilde{x}),
\end{align*}
using the associativity in the 2nd step and $B$-bilinearity and unitality in the 3rd step. 

\begin{rem}
As mentioned, one could consider more general bisections by requiring left $B$-linearity  
up to an automorphism of $B$ and substitute the above condition (2) by 
$\sigma\big(s(b)t(b') (x \ot \tilde{x}) \big)= \phi(b) \sigma(x \ot \tilde{x})b'$, and the condition (3) by
$\sigma\big((x \ot \tilde{x})(y \ot \tilde{y})\big)=\sigma\big((x \ot \tilde{x}) s(\phi^{-1}(\sigma(y \ot \tilde{y})))\big)$, for $\phi \in\Aut(B)$.
We limit ourselves to the smaller class of bisections as in Definition \ref{def:bisection}, 
that could be called vertical bisections, being dual to the vertical bisections on a classical gauge groupoid. 
Much of what follows can be adapted to the general bisections in a direct, if technically and notationally quite cumbersome, way. 
\end{rem}

The collection $\B(\C(A, H))$ of bisections of the bialgebroid $\C(A, H)$ is made a group by the convolution product of any two $\sigma_1$ bisections $\sigma_2$:    
\begin{align}\label{mulbis1}
\sigma_{1}\ast \sigma_{2}(x\ot  \tilde{x}) & := 
\sigma_{1}(\one{(x\ot \tilde{x})}) \, \sigma_{2}(\two{(x\ot \tilde{x})}) \nn \\ 
& \: = \sigma_{1}(\zero{x}\ot  \tuno{\one{x}}) \, \sigma_{2}(\tdue{\one{x}}\ot \tilde{x})
\end{align}
for any element $x\ot  \tilde{x}\in \C(A, H)$, recalling the $B$-coring coproduct \eqref{copro}. 

The product is well defined over the $B$-balanced tensor product since the bisections are $B$-bilinear; 
by the same reason $\sigma_{1}\ast\sigma_{2}$ is $B$-bilinear.  
We are left to show the associativity property (3) in the Definition \ref{def:bisection}. 
For simplicity write $X=x\ot  \tilde{x}$ and $Y=y\ot  \tilde{y}$ for elements of $\C(A, H)$ and 
use a Sweedler-like notation for the $B$-coring coproduct \eqref{copro}, 
$\Delta(X)=\one{X} \ot_B \two{X}$ and $\Delta(Y)=\one{Y} \ot_B \two{Y}$. Then, 
\begin{align*}
    \sigma_{1}\ast\sigma_{2}(XY)&=\sigma_{1}(\one{X}\one{Y})\sigma_{2}(\two{X}\two{Y})\\
    &=\sigma_{1}(\one{X}t(\sigma_{1}(\one{Y})))\sigma_{2}(\two{X}s(\sigma_{2}(\two{Y})))\\
    &=\sigma_{1}(\one{X}t(\sigma_{1}(\one{Y})\sigma_{2}(\two{Y})))\sigma_{2}(\two{X})\\
    &=\sigma_{1}(\one{X}s(\sigma_{1}(\one{Y})\sigma_{2}(\two{Y})))\sigma_{2}(\two{X})\\
    &=\sigma_{1}\ast\sigma_{2}(Xs(\sigma_{1}\ast\sigma_{2}(Y))),
\end{align*}
where the 2nd and 4th steps use \eqref{equ.bis}, and the 3rd step uses \eqref{eq:Tak.prod}. The product is associative since $\C(A, H)$ is coassociative as a $B$-coring.

The counit of the bialgebroid is a bisection by definition and one checks that 
$\varepsilon\ast\sigma=\sigma=\sigma\ast\varepsilon$ for any bisection $\sigma$ and $\varepsilon$ is the unit element. The inverse of the bisection $\sigma$ is
\beq\label{invobis1}
\sigma^{-1}(x\ot \tilde{x}) = x \, \sigma(\zero{\tilde{x}}\ot \tuno{\one{\tilde{x}}})\, \tdue{\one{\tilde{x}}},
\eeq
for $x \ot \tilde{x} \in \C(A, H)$ as we shall see below.
Indeed all group properties of $\B(\C(A, H))$ with the product \eqref{mulbis1} follow from the following proposition 
that parallels Proposition \ref{prop:gsogg}.

\begin{prop}\label{atob}
Let $B=A^{coH}\subseteq A$ be a faithfully flat Hopf--Galois extension,
and let $\C(A, H)$ be the corresponding Ehresmann--Schauenburg bialgebroid. 
There is a  group isomorphism $\alpha:\Aut_H(A)\to \B(\C(A, H))$ between gauge transformations and bisections.
\end{prop}
\begin{proof}
Firstly, given a bisection $\sigma \in \B(\C(A, H))$ we define a map  $F_{\sigma}: A\to A$ by
\beq\label{btog}
F_{\sigma}(a):=\sigma(\zero{a}\ot\tuno{\one{a}})\, \tdue{\one{a}} \, ,
\eeq
for any $a\in A$. 
This is well defined since $\sigma$ is right $B$-linear and 
$\zero{a}\otimes \tuno{\one{a}}\otimes_{B}\tdue{\one{a}}\in \C(A, H)\otimes_{B}A$. 
Clearly $F_{\sigma}|_{B}=\id_{B}$, and $F_{\sigma}$ is an algebra map:
\begin{align*}
    F_{\sigma}(aa')&=\sigma(\zero{a}\zero{a'}\otimes \tuno{\one{a'}}\tuno{\one{a}})\tdue{\one{a}}\tdue{\one{a'}}\\
    &=\sigma(\zero{a}\otimes\sigma(\zero{a}\otimes \tuno{\one{a'}})\tuno{\one{a}})\tdue{\one{a}}\tdue{\one{a'}}\\
    &=\sigma(\zero{a}\otimes \tuno{\one{a}})\tdue{\one{a}}\sigma(\zero{a'}\otimes \tuno{\one{a'}})\tdue{\one{a'}}\\
    &=F_{\sigma}(a)F_{\sigma}(a'),
\end{align*}
using \eqref{p2} for the 1st step, \eqref{equ.bis} for the 2nd, \eqref{p8} for the third. Also, $F_{\sigma}$ is $H$-equivariant:
\begin{align*}
\zero{F_{\sigma}(a)}\ot\one{F_{\sigma}(a)}&=\sigma(\zero{a}\ot\tuno{\one{a}})\zero{\tdue{\one{a}}}\ot\one{\tdue{\one{a}}}\\
&=\sigma(\zero{a}\ot\tuno{\one{a}})\tdue{\one{a}}\ot\two{a}\\
&=F_{\sigma}(\zero{a})\ot\one{a},
\end{align*}
where the 2nd step uses \eqref{p4}. Thus $F_{\sigma}\in\Aut_H(A)$.
 
\noindent
Conversely, let $F\in\Aut_H(A)$ a gauge transformation and define $\sigma_{F}\in \B(\C(A, H))$ by
\beq\label{gtob}
\sigma_{F}(a\ot \tilde{a}):=F(a)\tilde{a},
\eeq
for any $a\ot\tilde{a}\in \C(A, H)$. This is well defined since
\begin{align*}
\delta^A(F(a)\tilde{a}) &= \zero{(F(a)\tilde{a})}\ot \one{(F(a)\tilde{a})} =  \zero{F(a)}\zero{\tilde{a}}\ot \one{F(a)}\one{\tilde{a}} \\ & =F(\zero{a})\zero{\tilde{a}}\ot\one{a}\one{\tilde{a}} = F(a)\tilde{a}\ot 1_{H},
\end{align*}
using \eqref{ec2}. Clearly, $\sigma_{F}$ is unital and $B$-bilinear. Also, 
for any $a\otimes \tilde{a}, a'\otimes\tilde{a}'\in \C(A, H)$, 
\begin{align*}
    \sigma_{F}((a\otimes \tilde{a})(a'\otimes\tilde{a}'))&=F(a a')\tilde{a}' \tilde{a}=F(a)F(a')\tilde{a}' \tilde{a}\\
    &=F(aF(a')\tilde{a}')\tilde{a}=\sigma_{F}(a \sigma_{F}(a' \otimes\tilde{a}')\otimes\tilde{a}),
\end{align*}
where the 3rd step uses the fact that $F(a')\tilde{a}' \in B$. Thus $\sigma_{F}$ is a bisection.

\noindent
The correspondence is bijective: one easily checks that $\sigma_{F_{\sigma}} = \sigma $ and
$F_{\sigma_{F}}=F$. Also $\sigma_{\id_{A}} = \varepsilon$ 
and for any $a\ot \tilde{a}\in \C(A, H)$ we have,
\begin{align*}
\sigma_{G}\ast\sigma_{F}(a\ot \tilde{a})&=\sigma_{G}(\zero{a}\ot \tuno{\one{a}})  \,
\sigma_{F}(\tdue{\one{a}}\ot \tilde{a}) \\
&=G(\zero{a})\tuno{\one{a}}F(\tdue{\one{a}})\tilde{a}  =F\big(G(\zero{a})\tuno{\one{a}}\tdue{\one{a}} \big) \, \tilde{a} \\ & =F(G(a)) \, \tilde{a}=\sigma_{(G \circ F)}(a\ot \tilde{a}),
\end{align*}
where the 3th step uses that $G(\zero{a})\tuno{\one{a}}\otimes_{B}\tdue{\one{a}} = 
1 \ot_B G(\zero{a})\tuno{\one{a}} \tdue{\one{a}}$ similarly to \eqref{chi-eq}, and the 4th step uses \eqref{p5}.

Thus the correspondence is a group isomorphism. 
Via this we can get the inverse \eqref{invobis1} directly from \eqref{inv-F2}: 
$ \sigma^{-1}(a\otimes\tilde{a})=\sigma_{F_{\sigma}^{-1}}(a\otimes\tilde{a})=F_{\sigma}^{-1}(a)\tilde{a}=\zero{a}F_{\sigma}(\tuno{\one{a}})\tdue{\one{a}}\tilde{a} =aF_{\sigma}(\tilde{a})=a\sigma(\zero{\tilde{a}}\otimes\tuno{\one{\tilde{a}}})\tdue{\one{\tilde{a}}}$, 
where the 4th step uses \eqref{ec1}.
\end{proof}

We have already mentioned that gauge transformations for a noncommutative principal
bundles could be defined without them being algebra maps \cite{brz-tr}.
Mainly for the sake of completeness we record here a version of these via bialgebroids and bisections.
To distinguish them from the analogous concepts introduced previously, and for lack of a better name, we call them extended gauge transformation and extended bisections.

Thus following \cite{brz-tr}, 
the \textit{extended gauge group} $\Aut^{ext}_H(A)$ 
of a Hopf--Galois extension $B=A^{coH}\subseteq A$ 
consists of invertible $H$-comodule unital maps $F: A \to A$ which are such that 
$F(ba)=b F(a)$ for any $b\in B$ and $a\in A$. The group structure is map composition.

In parallel with this we have then the  following.
\begin{defi}\label{def:gebis}
Let $\C(A, H)$ be the Ehresmann--Schauenburg bialgebroid of the Hopf--Galois extension
$B=A^{co \, H}\subseteq A$. An \textit{extended bisection} is a convolution invertible, for the product \eqref{mulbis1}, unital $B$-bilinear 
map $\sigma: \C(A, H)\to B$. 
\end{defi}
Since extended bisections are $B$-bilinear, the product \eqref{mulbis1} is indeed well defined, and the counit $\varepsilon$ is still the unit element for the product. Thus the collection $\B^{ext}(\C(A, H))$ of extended bisections form a group, of which $\B(\C(A, H))$ is a subgroup. 
Notice that now \eqref{invobis1} is not the inverse in $\B^{ext}(\C(A, H))$ 
for the product in \eqref{mulbis1} since for an extended bisection we are not asking it to be a left character for the $B$-ring $(\C(A, H), s)$, and thus the expression in \eqref{invobis1} is not right $B$-linear. 

Finally, in analogy with Proposition \ref{atob} we have the following.
\begin{prop}\label{lem:gatb}
Let $B=A^{coH}\subseteq A$ be a faithfully flat Hopf--Galois extension with Ehresmann--Schauenburg bialgebroid $\C(A, H)$. Then there is a group isomorphism between extended gauge transformations $\Aut^{ext}_H(A)$ and extended bisections $\B^{ext}(\C(A, H))$.
\end{prop}
\begin{proof}
This uses the same methods as Proposition \ref{atob}.
Given $F\in\Aut^{ext}_H(A)$, define its image as in \eqref{gtob}: $\sigma_{F}(a\ot \tilde{a})=F(a)\tilde{a}$.
For all $b\in B$ we have,
\begin{align*}
\sigma_{F}( a\ot \tilde{a} b) = F(a) \tilde{a} \, b = \sigma_{F}(a\ot \tilde{a}) \, b \, , \quad 
    \sigma_{F}(b a \ot \tilde{a}) =F(b a)\tilde{a} = b F(a)\tilde{a} = b \sigma_{F}(a\ot \tilde{a}). 
\end{align*}   
Conversely, given $\sigma\in \B^{ext}(\C(A, H))$ define its image in $\Aut^{ext}_H(A)$ as in \eqref{btog}:
$F_{\sigma}(a)=\sigma(\zero{a}\ot\tuno{\one{a}})\tdue{\one{a}}$.
Then $F_{\sigma}(ba)=\sigma(b\zero{a}\ot\tuno{\one{a}})\tdue{\one{a}}=
b F_{\sigma}(a)$. The rest of the proof goes as that of Proposition \ref{atob}.
\end{proof}

\section{Hopf algebroids} \label{s5}
As examples we construct the Ehresmann--Schauenburg bialgebroid for the $\U(1)$ principal bundle over 
the quantum sphere in nonzero{ncu1b} and for the commutative not faithfully flat Hopf--Galois extension out of $\O(\mathrm{SL}(2))$ considered in Sect.~\ref{gtff}.
In both cases there is a suitable invertible antipode satisfying conditions \ref{hopbroid1} and \ref{hopbroid2} 
for a Hopf algebroid.
It is worth stressing that the results for $\O(\mathrm{SL}(2))$ in Sect.~\ref{gtff} and Sect.~\ref{biff} below show that at least for an  algebra of coinvariants which is commutative, the Hopf--Galois extension needs not be faithfully flat for all  constructions to be well defined. An analysis of the role of the faithful flatness in a general context will be reported elsewhere.

\subsection{The monopole bundle over the quantum sphere}
With reference to Sect.~\ref{ncu1b}, let us denote $H=\O(\U(1))$ coacting as in \eqref{rco} on the 
generators in \eqref{reltot} of $A=\ASL$ 
with algebra of coinvariants $B=\O(\sq)$ with generators in \ref{relbas}. 
From its definition $\C(A, H)=(A \ot A)^{coH}$, the bialgebroid is generated by elements
\beq\label{genu1b}
\alpha = a \ot d , \quad \gamma = - q^{-1} c \ot b,  \quad \widetilde{\alpha} = - q^{-1} a \ot b,  \quad \widetilde{\gamma} = c \ot d 
\eeq
and their `conjugated': 
\beq\label{genu1bb}
\delta   = d \ot a , \quad \beta   = - q^{-1} b \ot c,  \quad \widetilde{\beta} = d \ot c , \quad 
\widetilde{\delta} = - q^{-1} b \ot a .  
\eeq
By using the expression \eqref{trau1} for the translation map and the relations \eqref{reltot}, one checks  
that the above generators $h \ot k $ satisfy the condition $\zero{h} \ot \tau(\one{h}) k = h \ot k \ot _B 1_A$ as it should be from the alternative description \eqref{ec1} of the bialgebroid $\C(A, H)$. For instance:
\begin{align*}
\alpha \quad \mapsto \quad a \ot \tau(z) d & = a \ot (d \ot_B a d - q^{-1} b \ot_B c d ) \nn \\
& = a \ot ( d a d - q^{-1} b  c d ) \ot_B 1_A  \nn \\ 
& = a \ot d ( a d - q b c ) \ot_B 1_A = a \ot d \ot_B 1_A = \alpha \ot_B 1_A . 
\end{align*}
\begin{align*}
\delta \quad \mapsto \quad d \ot \tau(z^{-1}) a & = d \ot (a \ot_B d a - q c \ot_B b a ) \nn \\ 
& = d \ot ( a d a - q b c a ) \ot_B 1_A  \nn \\
& = d \ot a (d a - q^{-1} b c ) \ot_B 1_A = d \ot a \ot_B 1_A = \delta \ot_B 1_A .
\end{align*}
Similar computations work for the other generators.  
A direct computation leads to
\begin{align*}
 & \beta \gamma + \widetilde{\delta} \widetilde{\gamma} = B_0 \ot 1 = s(B_0),  
 \qquad
 \beta \gamma + \widetilde{\beta} \widetilde{\alpha} = 1 \ot B_0 = t(B_0) , \nn\\
 & \alpha \widetilde{\delta} + q^2 \widetilde{\alpha} \beta = B_- \ot 1 = s(B_-) ,
 \qquad 
 \widetilde{\alpha} \delta  + q^2 \gamma  \widetilde{\delta} = 1 \ot B_- = t(B_-) , \nn \\
& \widetilde{\gamma} \delta  + q^2 \gamma \widetilde{\beta} = B_+ \ot 1 = s(B_+) , 
\qquad 
\alpha \widetilde{\beta} + q^2 \widetilde{\gamma} \beta = 1 \ot B_+ = t(B_+) ,
\end{align*}
with source $s : B \to \C(A, H)$ and target $t : B \to \C(A, H)$ maps respectively.

The eight generators in \eqref{genu1b} and \eqref{genu1bb} are not independent. Indeed, define 
\begin{align}\label{ABCD}
A = a \ot d - q \, b \ot c, \qquad B = b \ot a - q^{-1} a \ot b , \nn \\ 
C = c \ot d - q \, d \ot c, \qquad D = d \ot a - q^{-1} c \ot b .
\end{align} 
Then, a direct computation shows that
\beq\label{rel1}
D A -q^{-1} C B = 1 \ot 1 = A D -q B C. 
\eeq
Also, the sphere relations in \eqref{relbas} translate into 
 \beq\label{rel2-3}
 \big( B_{+}\,B_{-} - B_{0} \, ( 1 - B_{0} ) \big) \ot 1 = 0 = 1 \ot \big( B_{+}\,B_{-} - B_{0} \, ( 1 - B_{0} ) \big).
 \eeq
An alternative way to show the relations among the generators is to observe 
by a direct computation that there are four `circle' relations:
 \begin{align}
 \delta \alpha  &= (1-B_0)\ot (1-B_0) , \qquad \beta \gamma = B_0 \ot B_0, \nn \\
 \widetilde{\beta} \widetilde{\alpha} &= B_0 \ot (1-B_0) , \qquad \widetilde{\delta} \widetilde{\gamma} = (1-B_0) \ot B_0. 
 \end{align}
These in turn imply 
 \beq\label{sph-rel-bi1}
\delta   \alpha + \beta   \gamma + \widetilde{\beta} \widetilde{\alpha}  + \widetilde{\delta} \widetilde{\gamma} = 1 \ot 1 .
\eeq
 Notice that the above relations which survive the classical limit $q=1$, are constraints among the generators and not commutation relations. For the latter one has the following. 
\begin{lem}
For the product and structure as in Definition \ref{def:reb}, the generators in \eqref{genu1b} and \eqref{genu1bb} of the bialgebroid $\C(A, H)$ satisfy the relations:
\begin{align*}
\alpha \gamma &= q^2 \gamma \alpha , \qquad 
\alpha \widetilde{\alpha}  = q \widetilde{\alpha}  \alpha , \qquad 
\alpha \widetilde{\gamma} = q \widetilde{\gamma} \alpha  \nn \\ 
\alpha \beta   &= q^2 \beta \alpha , \qquad 
\alpha \widetilde{\beta} = q \widetilde{\beta} \alpha + (1-q^2) \widetilde{\gamma} \beta , \qquad 
\alpha \widetilde{\delta} = q \widetilde{\delta} \alpha + (1-q^2) \widetilde{\alpha} \beta \nn \\ 
\gamma \widetilde{\alpha}  & = q^{-1} \widetilde{\alpha}  \gamma , \qquad 
\gamma \widetilde{\gamma} = q^{-1} \widetilde{\gamma} \gamma , \qquad 
\gamma \widetilde{\beta} = q \widetilde{\beta} \gamma , \qquad 
\gamma \widetilde{\delta} = q \widetilde{\delta} \gamma ,  \nn \\
\widetilde{\alpha}  \widetilde{\gamma} &= \widetilde{\gamma} \widetilde{\alpha} , \qquad 
\widetilde{\alpha}  \widetilde{\delta} = q^2 \widetilde{\delta} \widetilde{\alpha}   ,
\end{align*}
as well as
\begin{align*}
\gamma \beta   & = \beta   \gamma , \qquad \widetilde{\gamma} \widetilde{\delta} = \widetilde{\delta} \widetilde{\gamma} + (1-q^2) \beta   \gamma , \qquad \widetilde{\alpha}  \widetilde{\beta} = \widetilde{\beta} \widetilde{\alpha}  + (1-q^2) \beta   \gamma \nn \\
\alpha \delta   & = \delta   \alpha  + (1-q^2) (\widetilde{\delta} \widetilde{\gamma} + \widetilde{\beta} \widetilde{\alpha} ) + (1-q^2)^2 \beta   \gamma 
\end{align*}
(There are also `conjugated' relations directly derived from the previous ones.)
As a consequence $BC = CB$, for the generators in \eqref{ABCD}.
\end{lem}
\begin{proof}
This is the result of a direct computation given the relations in \eqref{reltot}. 
\end{proof}

\begin{lem}
The bialgebroid $\C(A, H)$ has a structure of a Hopf algebroid with coproduct \eqref{copro} which results into:
\begin{align}\label{cop1ex}
\Delta(\alpha) &= \alpha \ot_B \alpha + \widetilde{\alpha} \ot_B  \widetilde{\gamma}, \qquad
\Delta(\widetilde{\alpha}) = \alpha \ot_B \widetilde{\alpha} + \widetilde{\alpha} \ot_B \gamma , \nn \\
\Delta(\gamma) &= \widetilde{\gamma} \ot_B \widetilde{\alpha} + \gamma \ot_B  \gamma, \qquad
\Delta(\widetilde{\gamma}) = \widetilde{\gamma} \ot_B \alpha + \gamma \ot_B \widetilde{\gamma} , \nn \\
\Delta(\delta) &= \delta \ot_B \delta +q^2 \widetilde{\beta} \ot_B  \widetilde{\delta}, \qquad
\Delta(\widetilde{\beta}) = \delta \ot_B \widetilde{\beta} + q^2 \widetilde{\beta} \ot_B \beta , \nn \\
\Delta(\beta) &= \widetilde{\delta} \ot_B \widetilde{\beta} + q^2 \beta \ot_B  \beta, \qquad
\Delta(\widetilde{\delta}) = \widetilde{\delta} \ot_B \delta + q^2 \beta \ot_B \widetilde{\delta} ;
\end{align}
counit \eqref{counit} which results into:
\begin{align}\label{cou1ex}
\varepsilon(\alpha) &= 1 - q^2 B_0, \qquad \varepsilon(\gamma) = B_0,  \qquad \varepsilon(\widetilde{\alpha}) =  B_- ,  \qquad \varepsilon(\widetilde{\gamma}) = B_+ \nn \\ 
\varepsilon(\delta) &= 1 - B_0 , \qquad \varepsilon(\beta) = B_0,  \qquad \varepsilon(\widetilde{\beta}) = q^{-1} B_+,  \qquad \varepsilon(\widetilde{\delta}) = q^{-1} B_- ;
\end{align}
and antipode $S=S^{-1}$: 
\beq\label{anti1ex}
S(\alpha) = \delta, \qquad S(\gamma) = \beta 
\qquad S(\widetilde{\alpha}) = \widetilde{\delta}, \qquad S(\widetilde{\gamma}) = \widetilde{\beta} .
\eeq
\end{lem}
\begin{proof}
The expressions for the counit are clear. For the coproduct, from Definition \ref{copro}, 
\begin{align*}
\Delta(\alpha) &= a \ot \tau(z) \ot d = a \ot (d \ot_B a \ot d - q^{-1} b \ot_B c \ot d ) \\
& = a \ot d \ot_B a \ot d - q^{-1} a \ot b \ot_B c \ot d \\ &
= \alpha \ot_B \alpha + \widetilde{\alpha} \ot_B  \widetilde{\gamma},
\end{align*}
using the expression \eqref{trau1} for the translation map. Similarly,
\begin{align*}
\Delta(\delta) & = d \ot \tau(z^{-1}) \ot a = d \ot (a \ot_B d \ot a - q c \ot_B b \ot a ) \nn \\
&= d \ot a \ot_B d \ot a - q a \ot c \ot_B b \ot a \nn \\
&= \delta \ot_B \delta +q^2 \widetilde{\beta} \ot_B  \widetilde{\delta}. 
\end{align*}
Similar computations work for the other generators.  
Finally, 
the antipode in \eqref{anti1ex} clearly satisfies \eqref{hopbroid1}. Then, on the one hand, for condition \eqref{hopbroid2} 
for the generator $\alpha$:
\begin{align*}
\one{S(\one{ \alpha })}\two{\alpha }\ot_B \two{S(\one{\alpha})}
& = \one{S(\alpha)} \alpha \ot_B \two{S(\alpha)} 
+ \one{S(\widetilde{\alpha})} \widetilde{\gamma} \ot_B \two{S(\widetilde{\alpha})} \\
& = \one{\delta} \alpha \ot_B \two{\delta} + \one{\widetilde{\delta}} \widetilde{\gamma} \ot_B \two{\widetilde{\delta}} \\
& =\delta\alpha\ot_{B}\delta+ q^2 \widetilde{\beta} \alpha \ot_{B} \widetilde{\delta} 
+ \widetilde{\delta} \widetilde{\gamma} \ot_{B}\delta+ q^2 \beta \widetilde{\gamma} \ot_B \widetilde{\delta}  \\
& = (da - q^{-1} bc ) \ot da \ot_B \delta + q^2 (da - q^{-1} bc ) \ot dc \ot_B \widetilde{\delta}  \\
& = 1 \ot 1 \ot_{B} d (a \delta + q^2 c \widetilde{\delta} ) \\
& = 1 \ot 1 \ot_{B} d (a d - q cb ) \ot a \\
& = 1 \ot 1 \ot_{B} d \ot a \\
& = 1 \ot 1 \ot_{B}\delta = 1 \ot 1 \ot_{B}S(\alpha) . 
\end{align*}
using the relations $d a- q^{-1} bc = 1$ and  $a d - q cb = 1$. On the other hand, being $S^{-1} = S$, 
\begin{align*}
\one{S(\two{\alpha})} \ot_B \two{S(\two{ \alpha })}\one{\alpha }
& = \one{S(\alpha)} \ot_B \two{S(\alpha)}  \alpha
+ \one{S(\widetilde{\gamma})} \ot_B \two{S(\widetilde{\gamma})} \widetilde{\alpha} \\
& = \one{\delta} \ot_B \two{\delta}  \alpha 
+ \one{\widetilde{\beta}} \ot_B \two{\widetilde{\beta}} \widetilde{\alpha}\\
& =\delta \ot_{B}\delta \alpha+ q^2 \widetilde{\beta} \ot_{B} \widetilde{\delta} \alpha 
+ \delta \ot_{B} \widetilde{\beta} \widetilde{\alpha} + q^2 \widetilde{\beta} \ot_B \beta \widetilde{\alpha} \\
& = \delta \ot_B da \ot (da - q^{-1} bc ) - q \widetilde{\beta} \ot_B ba \ot (da - q^{-1} bc ) \\
& = (\delta d - q \widetilde{\beta} b) a \ot_B 1 \ot 1 \\
& = d \ot (ad - q cb) a \ot_B 1 \ot 1 \\
& = d \ot a \ot_B 1 \ot 1 \\& = \delta \ot_B 1 \ot 1= S(\alpha) \ot_B 1 \ot 1, 
\end{align*}
using again the relations $d a- q^{-1} bc = 1$ and  $a d - q cb = 1$.
Similar computations go for the other generators. This concludes the proof.
\end{proof}

Again, Proposition \ref{atob} determines the group of bisections $\B(\C(A, H))$ out of gauge transformations worked out in Sect.~\ref{ncu1b}. Both equations \eqref{btog} and \eqref{gtob} are well defined. In particular for the equation \eqref{btog} the map
$$
A \ni h \mapsto \zero{h}\ot \tuno{\one{h}}\ot_{B}\tdue{\one{h}}\in \C(A, H)\ot_B A
$$
is well defined. Using the expression \eqref{trau1} for the translation map one gets on generators: 
\begin{align*}
& a \quad \mapsto \quad a \ot d \ot_B a - q^{-1} a \ot b \ot_B c = \alpha \ot_B a + \widetilde{\alpha} \ot_B c \\
& c \quad \mapsto \quad c \ot d \ot_B a - q^{-1} c \ot b \ot_B c = \widetilde{\gamma} \ot_B a + \gamma \ot_B c \\
& b \quad \mapsto \quad b \ot a \ot_B d - q b \ot c \ot_B b = -q \widetilde{\delta} \ot_B d + q^2 \beta \ot_B b \\
& d \quad \mapsto \quad d \ot a \ot_B d - q d \ot c \ot_B b =  \delta \ot_B d - q \widetilde{\beta} \ot_B b .
\end{align*}
Given the generic gauge transformation in \eqref{qgt}, 
formula \eqref{gtob} determines a generic bisection. This is then given on generators by
\begin{align*}
\sigma(\alpha) & = X ad, \quad \sigma(\gamma) = -q^{-1} X cb, 
\quad \sigma(\widetilde{\alpha}) = -q^{-1} X ab, \quad\sigma(\widetilde{\gamma}) = X cd  \\
\sigma(\delta) &= X^{-1} da, \quad \sigma(\beta) = -q^{-1} X^{-1} bc, \quad \sigma(\widetilde{\beta}) = X^{-1} dc, \quad\sigma(\widetilde{\delta}) = -q^{-1} X^{-1} ba . 
\end{align*}
As before in Sect.~\ref{ncu1b}, $X$ is any map from the sphere $\mathrm{S}^2 \to \IC^*$ when $q=1$, while is any non vanishing element $X\in \IC^*$ when $q\not=1$.

%

\subsection{A commutative not faithfully flat example}\label{biff}
We work out the bialgebroid, in fact a Hopf algebroid, of the flat but not faithfully flat Hopf--Galois extension that we considered in Sect.~\ref{gtff}. Again there is an antipode leading to a Hopf algebroid.

Referring to that section, the Hopf algebra $H=\IC[x]$ with $x$ a primitive element, 
coacts as in \eqref{conff} on the algebra
$A = \O(\mathrm{SL}(2)) = \IC[a, b, c, d]/(ad - bc -1)$  with 
algebra of coinvariants $B= \IC[c, d]$. 
The bialgebroid $\C(A, H)=(A \ot A)^{coH}$ is generated by elements
\begin{align}\label{gen2ex}
\alpha &= a \ot d - c \ot b,  \qquad  \beta = b \ot d - d \ot b,  \nn \\
\gamma &= c \ot a - a \ot c,  \qquad \delta = d \ot a - b \ot c ,
\end{align}
easily seen to be invariant for the diagonal coaction.  The following is easily established.
\begin{lem}
For the product and structure as in Definition \ref{def:reb}, the generators above of the bialgebroid $\C(A, H)$ commute with each other while satisfying a sphere relation:
\beq\label{sph-rel-bi2}
\alpha \, \delta - \beta \gamma = 1 \ot 1 .
\eeq
Furthermore they also give
\begin{align}\label{abcdst}
\alpha c + \gamma d &= c \ot 1 = s(c) , \qquad \beta c + \delta d = d \ot 1 = s(d) , \nn \\
c \delta - d \gamma &= 1 \ot c = t(c) , \qquad d \alpha + - c \beta = 1 \ot d = t(d) , 
\end{align}
with source $s : B \to \C(A, H)$ and target maps $t : B \to \C(A, H)$.
\end{lem}

Also in the present case the bialgebroid $\C(A, H)$ can be given as in \eqref{ec1} since the generators 
$h \ot k $ satisfy the condition $\zero{h} \ot \tau(\one{h}) k = h \ot k \ot _B 1_A$. 
Using the expression \eqref{invnff} for the translation map, the coinvariance of $c$ and $d$ 
to pass them over the balanced tensor product, 
and the relation $ad - bc =1$, one computes for instance, 
\begin{align*}
\alpha = a \ot d - c \ot b \quad \mapsto \quad & \, \big(a \ot \tau(1) + c \ot \tau(z) \big) d- c \ot \tau(1) b   \\
& = a \ot d \ot_B 1 + c \ot a \ot_B bd - c \ot b \ot_B a d - c \ot 1 \ot_B b   \\
& = a \ot d \ot_B 1 + c \ot ad \ot_B b - c \ot b \ot_B a d - c \ot 1 \ot_B b   \\
& = a \ot d \ot_B 1 + c \ot bc \ot_B b - c \ot b \ot_B a d   \\
& = a \ot d \ot_B 1 + c \ot b \ot_B bc - c \ot b \ot_B a d   \\
& = (a \ot d  - c \ot b) \ot_B 1 = \alpha \ot_B 1 ,
\end{align*}
and similarly for the other generators. 

\begin{lem}
The bialgebroid $\C(A, H)$ has a structure of a Hopf algebroid with 
coproduct \eqref{copro} which results into:
\begin{align}\label{copr2ex}
\Delta(\alpha) &= \alpha \ot_B \alpha + \gamma \ot_B \beta, \qquad  
\Delta(\delta) = \delta \ot_B \delta + \beta \ot_B \gamma , \nn \\
\Delta(\beta) &= \beta \ot_B \alpha + \delta \ot_B \beta, \qquad  
\Delta(\gamma) = \gamma \ot_B \delta + \alpha \ot_B \gamma, 
\end{align}
counit \eqref{counit} which results into:
\beq\label{cou2ex}
\varepsilon(\alpha) = \varepsilon(\delta) = 1, \qquad \varepsilon(\beta) = \varepsilon(\gamma) = 0, 
\eeq
and antipode:
\beq\label{anti2ex}
S(\alpha) = \delta, \qquad S(\delta) = \alpha, \qquad S(\beta) = -\beta , \qquad S(\gamma) = -\gamma .
\eeq
\end{lem}
\begin{proof}
The form of the counit is clear. For the coproduct, from Definition \eqref{copro} and translation map \eqref{invnff}, 
one computes on the generators \eqref{gen2ex}, 
$$
\Delta(\alpha) = a \ot 1 \ot_B 1 \ot d + c \ot a \ot_B b \ot d - c \ot b \ot_B a \ot d - c \ot 1 \ot_B 1 \ot b 
$$
while, crossing $d$ and $c$ over the balanced tensor product,  
\begin{align*}
\alpha \ot_B \alpha & + \gamma \ot_B \beta \nn \\
&= a \ot 1 \ot_B d a \ot d - a \ot 1 \ot_B dc \ot b - c \ot b \ot_B a \ot d + c \ot b c \ot_B 1 \ot b \nn \\ 
& \quad + c \ot a \ot_B b \ot d - c \ot a d \ot_B 1 \ot b - a \ot 1 \ot_B c b \ot d + a \ot 1 \ot_B c d \ot b ,
\end{align*}
which coincides with the previous expression when using the relation $ad-bc=1$. Also, 
$$
\Delta(\beta) = b \ot 1 \ot_B 1 \ot d + d \ot a \ot_B b \ot d - d \ot b \ot_B a \ot d - d \ot 1 \ot_B 1 \ot b 
$$
which is the same as  
\begin{align*}
\beta \ot_B \alpha &+ \delta \ot_B \beta \nn \\
&= b \ot 1 \ot_B d a \ot d - b \ot 1 \ot_B dc \ot b - d \ot b \ot_B a \ot d + d \ot b c \ot_B 1 \ot b \nn \\ 
& \quad + d \ot a \ot_B b \ot d - d \ot a d \ot_B 1 \ot b - b \ot 1 \ot_B c b \ot d + b \ot 1 \ot_B c d \ot b .
\end{align*}
Similar computations work for the remaining generators. 
Finally, 
the antipode in \eqref{anti2ex} clearly satisfies \eqref{hopbroid1}. Then, on the one hand, for condition \eqref{hopbroid2}
for the generator $\alpha$:
\begin{align*}
\one{S(\one{ \alpha })}\two{\alpha }\ot_B \two{S(\one{\alpha})}
& = \one{S(\alpha)} \alpha \ot_B \two{S(\alpha)} + \one{S(\gamma)} \beta \ot_B \two{S(\gamma)} \\
& = \one{\delta} \alpha \ot_B \two{\delta} - \one{\gamma} \beta \ot_B \two{\gamma} \\
& =\delta\alpha\ot_{B}\delta+\beta\alpha\ot_{B}\gamma -\gamma\beta\ot_{B}\delta - \alpha\beta\ot_{B}\gamma  \\
& = 1 \ot 1 \ot_{B}\delta = 1 \ot 1 \ot_{B}S(\alpha),
\end{align*}
using the relation $\alpha \, \delta - \beta \gamma = 1 \ot 1$. On the other hand, being $S^{-1} = S$, 
\begin{align*}
\one{S(\two{\alpha})} \ot_B \two{S(\two{ \alpha })}\one{\alpha }
& = \one{S(\alpha)} \ot_B \two{S(\alpha)} \alpha + \one{S(\beta )} \ot_B \two{S(\beta )} \gamma \\
& = \one{\delta} \ot_B \two{\delta} \alpha - \one{\beta} \ot_B \two{\beta} \gamma \\
& = \delta \ot_{B}\delta \alpha +\beta \ot_{B}\gamma \alpha - \beta\ot_{B} \alpha\gamma - \delta \ot_{B} \beta \gamma \\
& = \delta \ot_{B} 1 \ot 1 = S(\alpha) \ot_{B} 1 \ot 1 .
\end{align*}
Similar computations go for the other generators. This concludes the proof.
\end{proof}

Finally, once again both equations \eqref{btog} and \eqref{gtob} are well defined since the map
$$
A \ni h \mapsto \zero{h}\ot \tuno{\one{h}}\ot_{B}\tdue{\one{h}}\in \C(A, H)\ot_B A
$$
is well defined. Indeed, using once more the expression \eqref{invnff} 
for the translation map one gets on the generators $a$ and $b$ 
(the not coinvariant ones): 
\begin{align}
a \quad \mapsto \quad & \, a \ot 1 \ot_B 1 + c \ot ( a \ot_B b - b \ot_B a ) \nn \\ 
& =  a \ot 1 \ot_B (a d - bc) + c \ot ( a \ot_B b - b \ot_B a ) \nn \\
& =  (a \ot d - c \ot b) \ot_B a + (c \ot a - a \ot c) \ot_B b \nn \\
& =  \alpha \ot_B a + \gamma \ot_B b , \label{ato}
\end{align}
having inserted $a d - bc = 1$ and using the coinvariance of $c$ and $d$ to cross them over the balanced tensor product.
Similarly for $b$ one finds:
\begin{align}
b \quad \mapsto \quad & b \ot 1 \ot_B 1 + d \ot ( a \ot_B b - b \ot_B a ) \nn \\
& = (d \ot a - b \ot c) \ot_B b + (b \ot d - d \ot b) \ot_B a \nn \\
& =  \beta \ot_B a + \delta \ot_B b . \label{bto}
\end{align}
Thus Proposition 
\ref{atob} relates gauge transformations of the Hopf--Galois extension to the group of bisections of the bialgebroid.
Let us illustrate this isomorphism by computing explicitly the identifications 
$F_{\sigma_F}=F$ and $\sigma_{F_\sigma} = \sigma$. 
Firstly, with definitions \eqref{btog} and \eqref{btog}, using \eqref{ato}, one computes
\begin{align}
F_{\sigma_F}(a) = \sigma_F(\alpha) a + \sigma_F(\gamma) b &= (F(a) d - c b) a + (c a - F(a) c) b \nn \\
& = F(a) (ad - bc) = F(a) .
\end{align}
A similar computation goes with the generator $b$, the statement being trivial 
for the coinvariant elements $c$ and $d$. Conversely, 
\begin{align}
\sigma_{F_\sigma}(\alpha) = F_\sigma(a) d - c b & = (\sigma(\alpha) a + \sigma(\gamma) b ) d - bc \nn \\
& = \sigma(\alpha) (1 + bc ) + \sigma(\gamma) b d - bc \nn \\
& = \sigma(\alpha) + (\sigma(\alpha c) + \sigma(\gamma) d) b - cb \nn \\
& = \sigma(\alpha) + (\sigma(\alpha c + \gamma d)) b - cb \nn \\
& = \sigma(\alpha) + \sigma(c \ot 1) b - cb = \sigma(\alpha) + c b - cb = \sigma(\alpha)
\end{align}
using the first relation in \eqref{abcdst}. Similar computations go with the remaining generators.
As an example, one gets for the gauge transformations in \eqref{gu1b} the bisection:
\beq
\sigma_F(\alpha) = 1 + h \, cd , \qquad \sigma_F(\beta) = h \, d^2 , 
\qquad \sigma_F(\gamma) = - h \, c^2 , \qquad \sigma_F(\delta) = 1 - h \, dc , 
\eeq
for $h$ an arbitrary element in $B$.

\section{Galois objects}\label{s6}
We shall now consider \emph{Galois objects} of a Hopf algebra $H$.
Such an object could be thought of as a noncommutative principal bundle over a point.
In contrast to the classical result that any fibre bundle over a point is trivial,
the set $\mbox{Gal}_H(\IC)$ of isomorphic classes of $H$-Galois objects needs
not be trivial (see  \cite{Bich}, \cite{kassel-review}). We shall illustrate later on this non-triviality with examples coming from group algebras and Taft algebras.

As already mentioned, the results of this section could be and have been obtained in an abstract and categorical way. Here we re-obtained them in a more explicit and more workable way, having in mind  
potential application to noncommutative  gauge theory.

\subsection{The bialgebroid of a Galois object}

\begin{defi}\label{Galoisobject}
Let $H$ be a Hopf algebra. A \textit{Galois object} of $H$ is an
$H$-Hopf--Galois extension $A$ of the ground field $\IC$.
\end{defi}
Thus, for a Galois object the coinvariant subalgebra is the ground field $\IC = A^{co \, H}$.
With coaction $\delta^A : A \to A\ot H$, $\delta^A(a) = \zero{a}\ot\one{a}$,
and translation map $\tau : H \to A \ot A$, $\tau(h)=\tuno{h}\ot\tdue{h}$,
for the Ehresmann--Schauenburg bialgebroid of a Galois object, being $B=\IC$, one has
(see also \cite[Def. 3.1]{schau}):
\begin{align*}
\C(A, H)
& = \{a\ot  \tilde{a}\in A\ot  A \, : \,\, \zero{a}\ot  \zero{\tilde{a}}\ot  \one{a}\one{\tilde{a}}=a\ot  \tilde{a}\ot  1_H \} \\ 
&=\{a\ot  \tilde{a}\in A\ot  A : \,\, \zero{a}\ot \tuno{\one{a}} \ot\tdue{\one{a}} \tilde{a} = a\ot \tilde{a} \ot 1_A\}  \, . 
\end{align*}
The coproduct \eqref{copro} and counit \eqref{counit} become $\Delta_\C(a\ot  \tilde{a}) = \zero{a} \ot \tuno{\one{a}} \ot \tdue{\one{a}} \ot \tilde{a}$,
and $\varepsilon_\C(a\ot  \tilde{a})=a\tilde{a} \in \IC$ respectively, for any $a\ot\tilde{a}\in \C(A, H)$. But now there is also an antipode \cite[Thm. 3.5]{schau} given, for any $a\ot\tilde{a}\in \C(A, H)$, by
\beq\label{anti}
S_{\C}(a\ot \tilde{a}):=\zero{\tilde{a}}\ot \tuno{\one{\tilde{a}}}a\tdue{\one{\tilde{a}}} \, .
\eeq
Thus the Ehresmann--Schauenburg bialgebroid of a Galois object is a Hopf algebra.

Recall that an $(A,H)$-relative Hopf module $M$ is a $H$-comodule with a compatible $A$-module structure.
That is the action is a morphism of $H$-comodules such that
$ \delta^M( ma) = \zero{m} \zero{a} \ot \one{m} \one{a}$ for all $a \in A$, $m\in M$. 
The multiplication induces an isomorphism \cite{HJSc90},  
$$
M^{co \, H} \ot A \to M,
$$
whose inverse is $M \ni m \mapsto \zero{m} \tuno{\one{m}} \ot \tdue{\one{m}} \in M^{co \, H} \ot A$ 
(see \cite[eq. 2.7]{schau2}).
Then, given that $\C(A, H)=(A\ot A)^{coH}$, this yields an isomorphism
$$
 A \ot A \simeq \C(A, H) \ot A \, , \qquad \widetilde{\chi}(a\ot  \tilde{a}) = \zero{a} \ot \tuno{\one{a}} \ot \tdue{\one{a}} \tilde{a} \, .
$$
 We finally collect some results of \cite{schau} (see  Lemma 3.2 and Lemma 3.3) in the following:
 \begin{lem}\label{universal}
 Let $H$ be a Hopf algebra, and $A$ a Galois object of $H$.
 There is a right $H$-equivariant algebra map $\delta^{\C} : A \to \C(A, H) \ot A$ given by
 $$
\delta^{\C}(a) = \zero{a} \ot \tuno{\one{a}} \ot \tdue{\one{a}}
$$
 which is universal in the following sense: Given an algebra $M$ and a $H$-equivariant algebra map
 $\phi : A \to M \ot A$, there is a unique algebra map $\Phi :  \C(A, H) \to M$ such that
 $\phi = (\Phi \ot \id_A) \circ \delta^{\C}$. Explicitly, $\Phi(a \ot  \tilde{a}) \ot 1_A = \phi(a) \tilde{a}$.
\end{lem}

Being algebra maps, now bisections are characters of the Hopf algebra
$\C(A, H)$ with product in \eqref{mulbis1} and inverse in \eqref{invobis1} that, 
with the antipode in \eqref{anti}
is written $\sigma^{-1} = \sigma \circ S_{\C}$, as it is the case for characters. From Proposition \ref{atob} we have the isomorphism
\beq\label{isogen}
\Aut_{H}(A) \simeq \B(\C(A, H)) = \mbox{Char}(\C(A, H)).
\eeq
This recover the result \cite[Cor. 3.1.4]{sch04}. 
If $H$ is the Hopf algebra of a compact quantum group $G$
this isomorphism is the result of \cite{AC20} that equivariant endomorphisms of $H$ 
are automorphisms and are `translations' by elements of the largest classical subgroup of $G$. 

For extended bisections and automorphisms Proposition \ref{lem:gatb} gives an isomorphism,
\beq\label{isogenext}
\Aut^{ext}_H(A) \simeq \B^{ext}(\C(A, H)) = \mbox{Char}^{ext}(\C(A, H)) \, ,
\eeq
with $\mbox{Char}^{ext}(\C(A, H))$ the group of convolution invertible unital maps
$ \phi: \C(A, H) \to \IC$.

\begin{exa}
Any Hopf algebra $H$ is a $H$-Galois object with its coproduct as coaction.
Then $H$ is isomorphic to the corresponding left bialgebroid $\C(H, H)$.

Indeed, 
if $H$ is a Hopf algebra with coproduct $\Delta(h)=\one{h}\ot\two{h}$, for the corresponding coinvariants: 
$\one{h}\ot\two{h}=h\ot 1$, we have $\varepsilon(\one{h})\ot\two{h}=\varepsilon(h)\ot 1$, this imply $h=\varepsilon(h)\in \IC$ and
$H^{co \, H}=\IC$. 
Moreover, the canonical Galois map $\chi :  g \ot h \mapsto  g \one{h} \ot \two{h}$ is bijective with inverse
given by  
$\chi^{-1}( g \ot h ):=  g \, S(\one{h})\ot\two{h}$. Thus $H$ is a $H$-Galois object.

With $A=H$, the corresponding left bialgebroid becomes
\begin{align*}
\C(H, H) & = \{ g \ot  h \in H \ot H : \,\, \one{g} \ot \one{h} \ot \two{g} \two{h} = g \ot h \ot  1_H \}  \\
& = \{g \ot h \in H\ot H : \,\, \one{g} \ot S(\two{g}) \ot \three{g} h = g \ot h \ot 1_A \} . 
\end{align*}
We have a linear map $\phi:\C(H, H) \to H$ given by $\phi(g\ot h):=g \, \varepsilon(h)$. 
The map $\phi$ has inverse $\phi^{-1}: H\to\C(H, H)$, defined by $\phi^{-1}(h):=\one{h}\ot S(\two{h})$. 
This is well defined since 
\begin{align*}
\Delta^{H \ot H}(\one{h}\ot S(\two{h}))=\one{h}\ot S(\four{h})\ot \two{h}S(\three{h})=\one{h}\ot S(\two{h})\ot 1_{H}, 
\end{align*}
showing that $\one{h}\ot S(\two{h}) \in\C(H, H)$. 
 Moreover,  
 $\phi(\phi^{-1}(h))=\phi(\one{h}\ot S(\two{h}))=h$, 
and
$\phi^{-1}(\phi(g\ot h)) = \varepsilon(h) \, \phi^{-1}(g) = \varepsilon(h) \, \one{g}\ot S(\two{g}) = g \ot h$.
The map $\phi$ is an algebra map: 
\begin{align*}
\phi((g \ot  h) \bullet_{\C} (g' \ot  h')) =
\phi (g g' \ot  h' h) = gg'  \varepsilon(h') \varepsilon(h) = \phi(g \ot  h) \bullet_{\C} \phi(g' \ot  h').
\end{align*}  
It is also a coalgebra map: 
\begin{align*}
(\phi\ot\phi)(\Delta_{\C}(g\ot h))&=(\phi\ot\phi)(\one{g}\ot \tuno{\two{g}}\ot \tdue{\two{g}}\ot h)\\
&=(\phi\ot\phi)(\one{g}\ot S(\two{g})\ot \three{g}\ot h)\\
&=\one{g}\ot\two{g} \, \varepsilon(h) \Delta_{H}(\phi(g\ot h));
\end{align*}
\begin{align*}
\varepsilon_{\C}(g\ot h) = gh = \varepsilon_{H}(gh)=\varepsilon_{H}(g)\varepsilon_{H}(h)=\varepsilon_{H}(\phi(g\ot h)) \, .
\end{align*}
\end{exa}

\subsection{Group Hopf algebras} \label{gradedalg} 
With a cocommutative Hopf algebra $H$ and $A$ a Galois object for $H$,
the bialgebroid $\C(A, H)$ is isomorphic to $H$ as a Hopf  algebra \cite[Rem. 3.8; Thm. 3.5]{schau}.
We work out some of the details for the case of a group algebra whose Galois objects 
are classified by group cohomology \cite[Chap. 8]{mont}.

Let $H=\IC[G]$ be a group algebra and let $A=\oplus_{g \in G} \, A_g$ be a strongly 
$G$-graded algebra. If $A$ is a $\IC[G]$-Galois object, that is $A_e=\IC$, each component $A_g$ 
is one-dimensional. If we pick a  element $u_g$ in each $A_g$, 
the multiplication of $A$ is determined by the products $u_g u_h$ for each pair $g, h$ of elements of $G$: 
\beq\label{mult}
u_{g}u_{h}=\lambda(g, h) \, u_{gh}
\eeq 
for a non vanishing $\lambda(g, h)\in \IC$. 
Associativity of the product requires that $\lambda$ satisfies a 2-cocycle condition, that is for any $g, h\in G$, 
$$
\lambda(g, h)\lambda(gh, k)=\lambda(h, k)\lambda(g, hk). 
$$
With a different nonzero element $v_g \in A_g$, we have $v_g = \mu(g) u_g$, 
for some nonzero $\mu(g) \in \IC$. The multiplication \eqref{mult} will become $v_{g}v_{h}=\lambda'(g, h)v_{gh}$ with 
$$
\lambda'(g, h) = \mu(g) \mu(h) (\mu(gh))^{-1} \lambda (g, h), 
$$
that is the two 2-cocycles $\lambda'$ and $\lambda$ are cohomologous. 
Thus the multiplication in $A$ depends only on the cohomology class of $\lambda \in H^2(G, \IC^{\times})$, the second cohomology group of $G$ with values in $\IC^{\times}$. Thus, equivalence classes of  $\IC[G]$-Galois objects are in bijective correspondence with the cohomology group $H^2(G, \IC^{\times})$.

\begin{exa}
For any cyclic group $G$ one has $H^2(G, \IC^{\times})=0$ and any corresponding 
$\IC[G]$-Galois object is trivial.
On the other hand, 
$H^2(\IZ^r, \IC^{\times}) = (\IC^{\times})^{r(r-1)/2}$
for the free abelian group of rank $r \geq 2$. 
Hence, there are infinitely many isomorphism classes of $\IC[\IZ^r]$-Galois objects (see \cite[Ex. 7.13]{kassel-review}).
\end{exa}

Being $H=\IC[G]$ cocommutative, as mentioned the bialgebroids $\C(A, H)$ are all 
isomorphic to $H$ as Hopf algebra. It is instructive to see this directly. For any $u_{g}\ot u_h \in \C(A, H)$
the coinvariance condition $u_{g} \ot u_h \ot g h = u_{g} \ot u_h \ot 1_H$, requires $h=g^{-1}$ so that $\C(A, H)$ is generated as vector space by elements $u_{g}\ot u_{g^{-1}}$, $g\in G$, with multiplication 
\beq\label{prog0}
(u_{g}\ot u_{g^{-1}})\bullet_{\C}(u_{h}\ot u_{h^{-1}}) = \lambda(g, h) \lambda(h^{-1}, g^{-1}) u_{gh}\ot u_{(gh)^{-1}}.
\eeq
Now, the cocycle 
$\Lambda(g, h) = \lambda(g, h) \lambda(h^{-1}, g^{-1})$ is trivial in $H^2(G, \IC^{\times})$. Indeed, 
$$
\Lambda(g, h) = \lambda(g, g^{-1}) \lambda(h, h^{-1}) / \lambda(gh, (gh)^{-1}) =  \mu(g) \mu(h) (\mu(gh))^{-1}, 
$$
with $\mu(g)=\lambda(g, g^{-1})$.
Consequently, by rescaling the generators $u_{g} \to v_g = \lambda(g, g^{-1})^{-\frac{1}{2}} \, u_{g}$
the multiplication rule \eqref{mult} becomes 
$v_{g}v_{h}=\lambda'(g, h) \, v_{gh}$,
with $\lambda'(g, h)= \Lambda(g, h)^{-\frac{1}{2}} \, \lambda(g, h)$ that we rename back to $\lambda(g, h)$.
As for the bialgebroid product in \eqref{prog0} one has,
$$
(v_{g}\ot v_{g^{-1}}) \bullet_{\C}(v_{h}\ot v_{h^{-1}}) = v_{gh}\ot v_{(gh)^{-1}},
$$
and the isomorphism $\Phi^{-1} : H \to \C(A, H) $ is simply $\Phi^{-1}(g)=\tau_L(g) = v_{g}\ot v_{g^{-1}}$.

The group of bisections $\B(\C(A, H))$ of $\C(A, H)$, and the gauge group $\Aut_H(A)$ of the Galois object $A$ coincide with the group of characters on $\IC[G]$, which is the same as $\Hom(G, \IC^{\times})$  the group (for point-wise multiplication) of group morphisms from  $G$ to $\IC^{\times}$.
Explicitly, since $F\in\Aut_H(A)$ is linear on $A$, on a basis  $\{v_{g}\}_{g\in G}$ of $A$, it is of the form
\begin{align*}
F(v_{g})=\sum\nolimits_{h\in G}f_{h}(g) v_{h},
\end{align*}
for complex numbers, $f_{h}(g) \in \IC$. Then, the $H$-equivariance of $F$, 
\begin{align*}
\zero{F(v_{g})}\ot\one{F(v_{g})}=F(\zero{v_{g}})\ot \one{v_{g}}=F(v_{g})\ot g,
\end{align*}
requires $F(v_{g})$ to be contained in to $A_{g}$ and we get $f_{h}(g) = 0$, if $h\not=g$ while $f_{g} := f_{g}(g) \in \IC^{\times}$ from the invertibility of $F$. Finally $F$ is an algebra map:
\begin{align*}
\lambda(g, h) f_{gh} v_{gh}=F(\lambda(g, h) \, v_{gh})=F(v_{g} v_{h}) 
= F(v_{g}) F(v_{h})=\lambda(g, h) f_{g} f_{h} \, v_{gh} \, ,
\end{align*}
implies $f_{gh}=f_{g}f_{h}$, for any $g, h\in G$. Thus we re-obtain that $\Aut_{H}(A)\simeq \Hom(G, \IC^{\times})$. 
Note that the requirement $F(v_e=1_A) = 1 = F_e$ implies that $F_{g^{-1}}=(F_{g})^{-1}$.

In contrast to this, the group $\Aut^{ext}_{H}(A)$ and then $\B^{ext}(\C(A, H))$ can be quite big.
If $F\in \Aut^{ext}_{H}(A)$,  that is one does not require $F$ to be an algebra map, the corresponding $f_g$ 
can take any value in $\IC^{\times}$ with the only condition that $f_e = 1$. 

\subsection{Taft algebras}\label{taftal}
Let $N\geq 2$ be an integer and let $q$ be a primitive $N$-th root of unity. 
The \emph{Taft algebra} $T_{N}$, introduced in \cite{Taft}, 
is a Hopf  algebra which is neither commutative nor cocommutative. Firstly, $T_{N}$
is the $N^{2}$-dimensional unital algebra generated 
by generators $x$, $g$ subject to the relations:
$$
x^{N}=0\, , \quad g^{N}=1\, , \quad xg - q\, gx=0 \, .
$$
It is a Hopf algebra with coproduct: 
\begin{align*}
\Delta(x):=1\ot x+x\ot g, \qquad \Delta(g):=g \ot g \, ;
\end{align*}
 counit: $\varepsilon(x):=0, \varepsilon(g):=1$, 
and antipode: $S(x): = - x g^{-1}, S(g) := g^{-1}$.
The four dimensional algebra $T_2$ is also known as the \emph{Sweedler algebra}.

For any $s\in \IC$, let $A_{s}$ be the unital algebra generated by elements $X, G$ with relations: 
$$
X^{N}=s\, , \quad G^{N} = 1\, , \quad X G - q \, GX=0\, . 
$$
The algebra $A_{s}$ is a right $T_{N}$-comodule algebra, with coaction defined by
\begin{align}\label{coTaft}
\delta^A(X):=1\ot x+ X \ot g,  \qquad \delta^A(G):=G \ot g.
\end{align}
The algebra of corresponding coinvariants is just the ground field $\IC$ and so $A_{s}$ is a $T_{N}$-Galois object.
It is known (see  \cite{Masuoka}, Prop. 2.17 and Prop. 2.22) that any $T_{N}$-Galois object is isomorphic to $A_s$ for some $s\in \IC$ and that any two such Galois objects $A_s$ and $A_t$ are isomorphic if and only if $s=t$. Thus the equivalence classes of  $T_{N}$-Galois objects are in bijective correspondence with the abelian group $\IC$. It is easy to see that the translation map of the coaction \eqref{coTaft} is given on generators by
\beq\label{trTaft}
\tau(g) = G^{-1} \ot G , \
\tau(x) = 1 \ot X - X G^{-1} \ot G .
\eeq
For the Ehresmann--Schauenburg bialgebroid $\C(A_s, T_N)$ we have then \cite[Cor. 2.4]{schau1}:
\begin{prop}\label{bitaft}
For any complex number $s$ there is a Hopf algebra isomorphism
$$
\Phi : \C(A_s, T_N) \simeq T_N .
$$
\end{prop}
\begin{proof}
Again we give a sketch of the explicit proof. It is easy to see that the elements 
$$
\Xi = X \ot G^{-1} - 1 \ot X G^{-1}  , \qquad  \Gamma = G \ot G^{-1}
$$
are coinvariants for the right diagonal coaction of $T_N$ on $A_s\ot A_s$ and that they generate 
$\C(A_s, T_N)=(A_s\ot A_s)^{co \, T_N}$ as an algebra. These elements satisfy the relations:
\begin{align*}
\Xi^N = 0, \quad \Gamma^N = 1, \quad \Xi \bullet_{\C} \Gamma = q \, \Xi \bullet_{\C} \Gamma  \, .
\end{align*}
The last two relations are easy to see. For the first one, one finds
\begin{align*}
\Xi^N & = X^N \ot G^{-N} + \sum_{r=1}^{N-1} c_r \, X^{N-r} \ot X^{r} G^{-N} + (-1)^N \ot (X G^{-1})^{N} \\
& = \big[ X^N \ot 1 + \sum_{r=1}^{N-1} c_r \, X^{N-r} \ot X^{r} + (-1)^N q^{\frac{n(n-1)}{2}}\ot X^{N} \big] \, G^{-N}
\end{align*}
when shifting powers of $G^{-1}$ to the right, for coefficients $c_r$ depending on $q$. Then, using the same methods as in \cite{Taft}, 
being $q$ a primitive $N$-th root of unity, all coefficients $c_r$ vanish and so
$\Xi^N = X^N \ot G^{-N} + (-1)^N \ot (X G^{-1})^{N}$ which then vanishes from $X^N=0$.

Thus $\Xi$ and $\Gamma$ generate a copy of the algebra $T_N$ and the isomorphism 
$\Phi$ maps $\Xi$ to $x$ and $\Gamma$ to $g$. 
The map $\Phi$ is a coalgebra map. Indeed,
$
\Delta(\Phi(\Gamma)) = \Delta(g) =g \ot g,
$
while,  
$
\Delta_\C(\Gamma) = \zero{G} \ot \tuno{\one{G}} \ot \tdue{\one{G}} \ot G^{-1} = G \ot G^{-1} \ot G \ot G^{-1} 
= \Gamma \ot \Gamma.
$
Thus $(\Phi \ot \Phi)(\Delta_\C(\Gamma)) = g \ot g = \Delta(\Phi(\Gamma))$. Similarly, 
$
\Delta(\Phi(\Xi)) = \Delta(x) = 1 \ot x + x \ot g,
$
while 
\begin{align*}
\Delta_\C(\Xi) & = \Delta_\C(X \ot G^{-1}) - \Delta_\C(1 \ot X G^{-1}) \\
& = \zero{X} \ot \tuno{\one{X}} \ot \tdue{\one{X}} \ot G^{-1} - 1 \ot 1 \ot 1 \ot X G^{-1} \\
& = 1 \ot \tuno{x} \ot \tdue{x} \ot G^{-1} + X \ot \tuno{g} \ot \tdue{g} \ot G^{-1}
- 1 \ot 1 \ot 1 \ot X G^{-1} \\
& = 1 \ot \Big(1 \ot X - X G^{-1} \ot G \Big) \ot G^{-1} + X \ot G^{-1} \ot G \ot G^{-1} - 1 \ot 1 \ot 1 \ot X G^{-1} \\
& = 1 \ot 1 \ot \Big( X \ot G^{-1} - 1 \ot X G^{-1} \Big) + \Big( X \ot G^{-1}  - 1 \ot X G^{-1} \Big) \ot G \ot G^{-1} \\
& = 1\ot \Xi + \Xi \ot \Gamma .
\end{align*}
  Thus $(\Phi \ot \Phi)(\Delta_\C(\Xi)) = 1 \ot x + x \ot g = \Delta(\Phi(\Xi))$. 
Finally: $\varepsilon_\C(\Gamma) = 1 = \varepsilon(g)$ and  $\varepsilon_\C(\Xi) = 0 = \varepsilon(x)$. 
This concludes the proof.
\end{proof}
The group of characters of the Taft algebra $T_N$ is the cyclic group $\IZ_N$: indeed any character $\phi$ must be such that  $\phi(x)=0$, while $\phi(g)^N=\phi(g^N)=\phi(1)=1$. Then for the group of gauge transformations of the Galois object $A_s$ --- the same as the group of bisections of the bialgebroid $\C(A_s, T_N)$ --- due to Proposition \ref{bitaft} we have,
$$
\Aut_{T_N}(A_s) \simeq \B(\C(A_s, T_N)) = \mbox{Char}(T_N) = \IZ_N. 
$$

An element $F$ of $\Aut^{ext}_{T_N}(A_s) \simeq \B^{ext}(\C(A_s, T_N)$, due to equivariance 
$\zero{F(a)}\ot\one{F(a)}=F(\zero{a})\ot \one{a}$ for any $a \in A_{s}$, can be given as a block diagonal matrix 
$$
F = \diag(M_1, M_2, \dots, M_{N-1}, M_N)
$$
with each $M_j$ a $N \times N$ invertible lower triangular matrix

$$
M_j = \begin{bmatrix}
1  & 0    & 0      &  \dots      & 0   & 0 \\
b_{21}  & a_{N-1} & 0 & \dots &   0  & 0 \\
b_{31}  & b_{32} & a_{N-2} & \ddots & \ddots &  \vdots  \\
  \vdots & \ddots & \ddots & \ddots & 0 & 0  \\
b_{N-1, 1} &  b_{N-1, 2}      & \ddots & \ddots & a_2 & 0 \\
b_{N 1} & b_{N 2}&  \dots      & b_{N, N-2}  & b_{N, N-1}      & a_1
\end{bmatrix}
$$

\noindent
All matrices $M_j$ have in common the diagonal elements $a_j$ (ciclic permuted) which are all different from zero for the invertibility of $M_j$. For the subgroup $\Aut_{T_N}(A_s)$ the $M_j$ are diagonal as well with 
$a_k=(a_1)^k$ and $(a_1)^N=1$ so that $M_j \in\IZ_N$. 
The reason all $M_j$ share the same diagonal elements (up to permutation) is the following: firstly, the `diagonal' form of the coaction of $G$ in \eqref{coTaft} implies that 
the image $F(G^k)$ is proportional to $G^k$, say $F(G^k)=\alpha_k G^k$ for some  constant $\alpha_k$. 
Then, due to the first term in the coaction of $X$ in \eqref{coTaft}, the `diagonal' component along the basis element $X^l G^k$ of the image $F(X^l G^k)$ is given again by $\alpha_k$ for any possible value of the index $l$.

Let us illustrate the construction for the cases of $N=2, 3$. Firstly, $F(1)=1$ since $F$ is unital. 
When $N=2$, on the basis $\{1, X, G, XG\}$, the equivariance $\zero{F(a)}\ot\one{F(a)}=F(\zero{a})\ot \one{a}$ for the coaction \eqref{coTaft} becomes
\begin{align*}
\zero{F(X)}\ot\one{F(X)}&=1 \ot x+F(X)\ot g,\\
\zero{F(G)}\ot \one{F(G)}&=F(G)\ot g,\\
\zero{F(XG)}\ot\one{F(XG)}&=F(G)\ot xg+F(XG)\ot 1.
\end{align*} 
Next, write $F(a)=f_{1}(a) + f_{2}(a)\, X+f_{3}(a)\, G + f_{4}(a)\, XG$, for complex numbers $f_{k}(a)$. 
And compute $\zero{F(a)}\ot \one{F(a)}=f_{1}(a)\, 1 \ot 1+ f_{2}(a)\, (1\ot x+ X\ot g)+f_{3}(a)\, G\ot g
+f_{4}(a)\, (G\ot xg+XG\ot 1)$. Then comparing generators, the equivariance gives
\begin{align*}
& f_1(X) = f_4(X) = 0 \\
& f_1(G) = f_2(G) = f_4(G)= 0 \\
& f_2(XG) = f_3(XG) = 0 ,
\end{align*}
while the remaining coefficients are related by the system of equations
\begin{align*}
f_{2}(X)\, (1\ot x+ X\ot g)+f_{3}(X)\, G\ot g &=1 \ot x+F(X)\ot g,\\
f_{3}(G)\, G\ot g &=F(G)\ot g,\\
f_{1}(XG)\, 1 \ot 1 + f_{4}(XG)\, (G\ot xg+XG\ot 1) &=F(G)\ot xg+F(XG)\ot 1.
\end{align*} 
One readily finds solutions
\begin{align*}
f_2(X) = 1 , \quad f_3(X) = \gamma , \quad f_1(XG) = \beta , \quad f_3(G) = f_4(XG) = \alpha 
\end{align*}
with $\alpha, \beta, \gamma$ arbitrary complex numbers. Thus a generic element $F$ of $\Aut^{ext}_{T_2}(A_s)$
can be represented by the matrix:
\begin{align}\label{autverT}
F : 
\begin{pmatrix}
1 \\ X G\\ G \\ X \\
\end{pmatrix} 
\quad \mapsto \quad
 \begin{pmatrix}
1 & 0 & 0 & 0\\
 \beta & \alpha & 0 & 0\\
0 & 0 & \alpha & 0\\
0 & 0 & \gamma & 1\\ 
\end{pmatrix} 
\begin{pmatrix}
1 \\ X G\\ G \\ X \\
\end{pmatrix} .
\end{align}
Asking $F$ to be invertible requires $\alpha \neq 0$.

On the other hand, any $F\in\Aut_{T_{2}}(A_{s})$ is an algebra map 
and so is determined by its values on the generators $G, X$.
From $F(G)= \alpha G$ and $F(X)= \gamma G + X$: 
requiring $s=F(X^2) = (\gamma G + X)^2 = \gamma + (GX + XG) + s$ yields $\gamma=0$; 
then $\beta + \alpha X G = F(XG) = \alpha XG$ yields $\beta=0$; and $1 = F(G^2) = (\alpha G)^2$ leads to 
re-obtain that $\Aut_{T_{2}}(A_{s})\simeq\mathbb{Z}_{2}$.

When $N=3$ a similar, if longer computation, gives for $\Aut^{exp}_{T_{3}}(A_{s})$ an eight parameter group with its elements of the following block diagonal form 
\begin{align*}
F :
\begin{pmatrix}
1 \\ X G^{2} \\ X^{2} G 
\end{pmatrix} 
\quad \mapsto \quad
 \begin{pmatrix}
1 & 0 & 0  \\
\beta & \alpha_2 & 0 \\
\eta & - q \delta & \alpha_1 \\ 
\end{pmatrix}  
\begin{pmatrix}
1 \\ X G^{2} \\ X^{2} G 
\end{pmatrix} .
\end{align*}
\begin{align*}
F :
\begin{pmatrix}
G \\ X \\ X^{2} G^2
\end{pmatrix} 
\quad \mapsto \quad
 \begin{pmatrix}
\alpha_1 & 0 & 0\\
\gamma & 1 & 0\\
\theta &- q\beta & \alpha_2 
\end{pmatrix}  
\begin{pmatrix}
G \\ X \\ X^{2} G^2
\end{pmatrix} .
\end{align*}
\begin{align*}
F :
\begin{pmatrix}
G^{2} \\ X G \\ X^{2} 
\end{pmatrix} 
\quad \mapsto \quad
 \begin{pmatrix}
\alpha_2 & 0 & 0 \\
\delta & \alpha_1 & 0 \\
\lambda & - q \gamma & 1 
\end{pmatrix}  
\begin{pmatrix}
G^{2} \\ X G \\ X^{2} 
\end{pmatrix} .
\end{align*}
One needs $\alpha_j \neq 0$, $j=1,2$ for invertibility. 

By going as before, for any $F\in\Aut_{T_{N}}(A_{s})$ one starts from it values on the generators, 
$F(G)= \alpha_1 G$ and $F(X)= \gamma G + X$, to conclude that $F$ is a diagonal matrix (in particular $F(X)=X$)
with $\alpha_2=(\alpha_1)^2$ and $1=(\alpha_1)^3$; thus $\Aut_{T_{3}}(A_{s})\simeq\mathbb{Z}_{3}$. 

\section{Crossed module structures on bialgebroids} \label{s7}

Automorphisms of a groupoid  with its natural transformations form a strict 2-group or, equivalently, a crossed module (see \cite[Def.~3.21]{ralf-chen}). The crossed module combines automorphisms and bisections. A bisection $\sigma$ is the 2-arrow from the identity morphism to an automorphism $Ad_\sigma$ and the composition of bisections can be viewed as the horizontal composition of 2-arrows. Then the  crossed module involves the product on bisections and the composition on automorphisms, and the group homomorphism from bisections to automorphisms together with the action of automorphisms on bisections by conjugation. 

In this section we quantise this construction for the Ehresmann--Schauenburg bialgebroid of a Hopf--Galois extension. We construct a crossed module for the bisections and the automorphisms of the bialgebroid, thus giving a generalization of a crossed module on a groupoid.
Notice that the antipode of the bialgebroid is not needed in the construction.

\subsection{Crossed modules and bisections} \label{cmhg}
A crossed module is the data $(M, N, \mu, \alpha)$ of two groups $M$, $N$ together with 
group homomorphisms $\mu: M\to N$ and $\alpha : N \to \Aut(M)$ such that, 
denoting $\alpha_{n} : M\to M$ for every $n\in N$, the following conditions are satisfied:
\begin{itemize}
    \item[(1)] $\mu(\alpha_{n}(m))=n\mu(m) n^{-1}$, \quad for any $n\in N$ and $m\in M$ ,
    ~\\
    \item[(2)] $\alpha_{\mu(m)}(m')=m \, m' \, m^{-1}$, \quad for any $m, m'\in M$.
\end{itemize}

We aim at proving the following.
\begin{thm}\label{thm. crossed module1}
Let $B=A^{coH}\subseteq A$ be a faithfully flat Hopf--Galois extension, with corresponding Ehresmann--Schauenburg bialgebroid $\C (A, H)$. Then there is a group morphism 
$Ad: \B(\C(A, H))\to\Aut(\C(A, H))$ and an action $\triangleright$ of $\Aut(\C(A, H))$ on $\B(\C(A, H))$ that give a crossed module structure to 
the pair $\big( \B(\C(A, H)), \Aut(\C(A, H)) \big)$. 
\end{thm}

We give the proof in a few lemmas.
\begin{lem}\label{ajaut}
With the hypothesis of Theorem \ref{thm. crossed module1}, consider a bisection $\sigma\in \B(\C(A, H))$
with $F_{\sigma} \in \Aut_H(A)$ the associated gauge transformation as in \eqref{btog}.  
Then, the map  
$Ad_\sigma: \C(A, H)\to \C(A, H)$ defined by
\begin{align} \label{Ad}
Ad_\sigma(a\ot\tilde{a}) & :=F_{\sigma}(a)\ot F_{\sigma}(\tilde{a}) \nn \\ & \:= 
\sigma(\zero{a}\ot\tuno{\one{a}})\tdue{\one{a}}\ot\sigma(\zero{\tilde{a}}\ot\tuno{\one{\tilde{a}}})\tdue{\one{\tilde{a}}}\, 
\end{align}
for any $a\ot\tilde{a} \in \C(A, H)$, is a vertical automorphism of $\C(A, H)$.
\end{lem}
\begin{proof}
Since $F_{\sigma}$ is an algebra automorphism so is $Ad_\sigma$. 
Then, for any $b\in B$ it is immediate to show that $Ad_\sigma(t(b)) = t (b)$ and 
$Ad_\sigma(s(b)) = s(b)$. So conditions \eqref{amoeba(i)} are satisfied.
The second one also shows that $Ad_\sigma$ is vertical, that is 
$\varepsilon \circ Ad_\sigma \circ s = \id_B$, and then
$\varepsilon\circ Ad_\sigma = \varepsilon$.
For the first condition in \eqref{amoeba(ii)} the
$H$-equivariance of $F_{\sigma}$ yields 
\beq\label{ad1}
(\Delta_{\C(A, H)}\circ Ad_\sigma) (a\ot\tilde{a}) 
= F_{\sigma}(\zero{a})\ot \tuno{\one{a}}\ot_{B}\tdue{\one{a}}\ot F_{\sigma}(\tilde{a}) 
\eeq
for any $a\ot\tilde{a}\in \C(A, H)$. On the other hand, 
\beq\label{ad2}
 (Ad_\sigma\ot_{B} Ad_\sigma)\circ \Delta_{\C(A, H)} (a\ot\tilde{a}) = F_{\sigma}(\zero{a})\ot F_{\sigma}(\tuno{\one{a}})\ot_{B}F_{\sigma}(\tdue{\one{a}})\ot F_{\sigma}(\tilde{a}).
\eeq
Now, for any $F \in \Aut_H(A)$,  given $h \in H$, one has 
$$
 F(\tuno{h}) \ot_{B}F(\tdue{h}) = \tuno{h}\ot_{B}\tdue{h} 
$$
as can be seen by applying the canonical map $\chi$ (an isomorphism) and using equivariance of $F$. 
Using it for the right hand sides of \eqref{ad1} and \eqref{ad2} shows that they coincide. Thus the left  
hand side expressions coincide and the first of \eqref{amoeba(ii)} is satisfied as well.
\end{proof}

\begin{rem}
The map $Ad_\sigma$ in \eqref{Ad} can also be written in the following useful ways:
\begin{align}\label{altad}
    Ad_\sigma(a\ot \tilde{a})& = \sigma(\zero{a}\ot\tuno{\one{a}})\tdue{\one{a}}\ot \tuno{\two{a}}\sigma^{-1}( \tdue{\two{a}} \ot \tilde{a}) \nn \\
   & = \sigma(\one{(a\ot \tilde{a})}) \triangleright \two{(a\ot \tilde{a})}  
\triangleleft \sigma^{-1} (\three{(a\ot \tilde{a})} ) 
 \end{align}
with the $B$-bimodule action \eqref{eq:rbgd.bimod}. 
Indeed, for $a\ot \tilde{a}\in \C(A, H)$, insert \eqref{p5} and compute:
\begin{align*}
    Ad_\sigma(a\ot \tilde{a})& = \sigma(\zero{a}\ot\tuno{\one{a}})\tdue{\one{a}}\ot \sigma(\zero{\tilde{a}}\ot \tuno{\one{\tilde{a}}})
   \tdue{\one{\tilde{a}}}
    \\        & = \sigma(\zero{a}\ot\tuno{\one{a}})\tdue{\one{a}}\ot \tuno{\two{a}}   \tdue{\two{a}} 
    \sigma(\zero{\tilde{a}}\ot  \tuno{\one{\tilde{a}}})\tdue{\one{\tilde{a}}} \\ 
    & = \sigma(\zero{a}\ot\tuno{\one{a}})\tdue{\one{a}}\ot \tuno{\two{a}} \sigma(\zero{\tilde{a}}\ot  \tuno{\one{\tilde{a}}})\tdue{\two{a}} \tdue{\one{\tilde{a}}}\\
    & = \sigma(\zero{a}\ot\tuno{\one{a}})\tdue{\one{a}}\ot \tuno{\two{a}} 
    \sigma^{-1}  
    ( \tdue{\two{a}} \ot \tilde{a}) \big)
      \\ & = \sigma(\one{(a\ot \tilde{a})}) \triangleright \two{(a\ot \tilde{a})} 
      \triangleleft \sigma^{-1} (\three{(a\ot \tilde{a})} ) .
\end{align*}
\end{rem}
One sees that $Ad_\sigma \circ Ad_\tau = Ad_{\tau \ast \sigma}$ for any
$\sigma$, $\tau \in \B(\C(A, H))$, while $(Ad_\sigma)^{-1} = Ad_{\sigma^{-1}}$ 
and $Ad_\varepsilon = \id_{\C(A, H)}$. Thus, 
$Ad$ is a group morphism $Ad: \B(\C(A, H))\to\Aut(\C(A, H))$.

Next, given an automorphism $(\Phi, \varphi)$ of $\C (A, H)$ 
with inverse $(\Phi^{-1}, \varphi^{-1})$, 
we define an action of $(\Phi, \varphi)$ on the group of bisections $\B(\C(A, H))$ as follows:
\beq\label{aaobs}
\Phi \triangleright \sigma := \varphi^{-1} \circ \sigma \circ \Phi 
\eeq
for any $\sigma\in \B(\C(A, H))$. 
The result is a well-defined bisection. It is clearly unital. Then, for $b, b' \in B$, $X \in \C(A, H)$, using \eqref{amoeba(i)} for the automorphism $\Phi$ and $B$-bilinearity of $\sigma$:
\begin{align*}
\Phi \triangleright \sigma \big(s(b)t(b') X \big) & = \varphi^{-1} \circ \sigma \circ \Phi \big(s(b)t(b') X \big) 
= \varphi^{-1} \circ \sigma \big( s(\varphi(b)) t(\varphi(b')) \Phi (X) \big) \\
& = \varphi^{-1} \big(\varphi(b) (\sigma \circ \Phi (X)) \varphi(b') \big) = b \, \big( \varphi^{-1} \circ \sigma \circ \Phi (X) \big) \, b'  \\
& = b \, (\Phi \triangleright \sigma (X))\, b' )
\end{align*}
and $\Phi \triangleright \sigma$ is $B$-bilinear as well. Finally, for all $X, Y \in \C(A, H)$, 
\begin{align*}
 \Phi \triangleright \sigma \big(X \, s( \Phi \triangleright \sigma (Y) \big) & = 
\varphi^{-1} \circ \sigma \circ \Phi \, \big(X \, s(\varphi^{-1} \circ \sigma \circ \Phi (Y)\big) =
\varphi^{-1} \circ \sigma \big( \Phi(X) \, s ( \sigma (\Phi (Y)) \big) \\
& = \varphi^{-1} \circ \sigma \big( \Phi(X) \Phi (Y) \big) = \varphi^{-1} \circ \sigma \circ \Phi (X Y) \\
& = \Phi \triangleright \sigma ( XY )
 \end{align*} 
using again \eqref{amoeba(i)} and the associativity of $\sigma$. Thus $\Phi \triangleright \sigma$ is associative.
One checks that
\beq\label{aaobs-1}
(\Phi\triangleright\sigma)^{-1} = \varphi^{-1} \circ \sigma^{-1} \circ \Phi .
\eeq

\begin{lem}\label{aaction}
Given an automorphism $(\Phi, \varphi)$ 
the action \eqref{aaobs} is in $\Aut(\C(A, H))$, the group automorphisms of $\B(\C(A, H))$.
\end{lem}
\begin{proof}
Let $\sigma, \tau\in \B(\C(A, H))$ and $X \in \C(A, H)$. We compute:
\begin{align*}
(\Phi\triangleright\tau)\ast(\Phi\triangleright\sigma)(X)
& =
(\Phi\triangleright\tau)(\one{X}) (\Phi\triangleright\sigma)(\two{X})\\
&=  \varphi^{-1} \circ \tau \circ \Phi (\one{X}) \, \, \varphi^{-1} \circ \sigma \circ \Phi(\two{X}) \\
&= \varphi^{-1} \big( \tau(\Phi(\one{X})) \, \sigma (\Phi(\two{X}) ) \big) = \varphi^{-1} \big(\tau(\one{\Phi(X)}) \, \sigma (\two{\Phi(X)}) \\
&= \varphi^{-1} \circ (\tau\ast\sigma)\circ \Phi(X) \big) \\
&= \Phi\triangleright(\tau\ast\sigma)(X)
\end{align*}
using the first equivariant condition \eqref{amoeba(ii)}. Also
\begin{align*}
\Phi \triangleright \varepsilon =  \varphi^{-1} \circ \varepsilon \circ \Phi =  \varphi^{-1} \circ \varphi \circ \varepsilon = \varepsilon.
\end{align*}
Finally, for any two automorphisms $(\Phi, \varphi)$ and $(\Psi, \psi)$ of $\C(A, H)$, we have
\begin{align*}
\Phi\triangleright(\Psi\triangleright(\sigma)) =   
\varphi^{-1} \circ \psi^{-1} \circ \sigma \circ \Psi \circ \Phi
= (\psi \circ \varphi)^{-1} \circ \sigma \circ \Psi\circ \Phi = (\Psi\circ \Phi)\triangleright\sigma.
\end{align*}
In particular $\Phi^{-1}\triangleright(\Phi\triangleright(\sigma))=\sigma$ and so the action is an automorphism 
of $\B(\C(A, H))$.
\end{proof}

\begin{lem}\label{d2}
For any automorphism $(\Phi, \varphi)$, and any $\sigma\in \B(\C(A, H))$ 
we have 
$$
Ad_{\Phi \triangleright \sigma}=\Phi^{-1} \circ Ad_\sigma \circ \Phi .
$$
\end{lem}
\begin{proof}
With $X \in \C(A, H)$, from \eqref{altad} we get
\begin{align}
    (Ad_\sigma \circ \Phi) (X) 
    = \sigma(\one{(\Phi (X))}) \triangleright \two{(\Phi (X))} 
   \triangleleft \sigma^{-1} \big(\three{(\Phi (X))}\big) , \label{11}
 \end{align}
while, using \eqref{aaobs} and \eqref{aaobs-1}, we have
\begin{align*}
    Ad_{\Phi\triangleright\sigma} (X) & = 
    \big( (\Phi\triangleright\sigma) (\one{X}) \big) \triangleright
    \two{X} \triangleleft \big( (\Phi\triangleright\sigma)^{-1} (\three{X}) \big) \\
    & = \big( \varphi^{-1}\circ \sigma (\Phi(\one{X}) ) \big) \triangleright \two{X} \triangleleft 
    \big(\varphi^{-1} \circ \sigma^{-1}( \Phi(\three{X})) \big) . 
 \end{align*}
Since $\Phi$ is a bimodule map: $\Phi(b  \triangleright X  \triangleleft \tilde{b}) = \varphi(b)  \triangleright \Phi(X)  \triangleleft \varphi(\tilde{b})$, 
for all $b,\tilde{b} \in B$, we get, 
\begin{align}
 (\Phi \circ Ad_{\Phi\triangleright\sigma}) (X)  
 = \sigma (\Phi(\one{X}) )  \triangleright \Phi(\two{X})  \triangleleft
\sigma^{-1} \big( \Phi(\three{X}) \big) .
  \label{33}
 \end{align}
The right hand sides of \eqref{11} and \eqref{33} are equal from the equivariance 
condition \eqref{amoeba(ii)}. 
\end{proof}

\begin{lem}\label{d3}
Let $\sigma$, $\tau\in \B(\C(A, H))$, then
$Ad_\tau \triangleright\sigma=\tau\ast\sigma\ast\tau^{-1}$.
\end{lem}
\begin{proof}
Recall that $Ad_\tau$ is vertical. With  $X \in \C(A, H)$, using definition \eqref{altad} we compute:
\begin{align*}
 Ad_\tau \triangleright\sigma(X) & = \sigma \circ Ad_\tau (X) ) = \sigma \big( \tau(\one{X})  
 \triangleright \two{X} \triangleleft \tau^{-1}  (\three{X}) \big)\\
& = \tau(\one{X})\, \sigma(\two{X}) \, \tau^{-1}(\three{X}) \\
& = \tau\ast\sigma\ast\tau^{-1}(X)
\end{align*}
where we used the definition \eqref{mulbis1} for 
the product. 
\end{proof}
Taken together the previous lemmas establish the content of Theorem \ref{thm. crossed module1} 
that is a crossed module structure for
$\big( \B(\C(A, H)), \Aut(\C(A, H)), Ad, \triangleright \big)$.

\begin{exa}\label{crossedmodule1}
Given a Hopf algebra $H$ and a character $\phi: H\to \mathbb{C}$, one defines a Hopf algebra automorphisms 
(see \cite[page 3807]{schau}) by
\begin{align}\label{autohopf}
\mbox{coinn}(\phi):H\to H , \quad \mbox{coinn}(\phi)(h):= \phi(\one{h}) \two{h} \phi(S(\three{h})),
\end{align}
for any $h\in H$. Recall that for a character $\phi^{-1} = \phi \circ S$.
The set $\mbox{CoInn}(H)$ of co-inner automorphisms of $H$ is a normal subgroup of the group
$\Aut_{\mbox{\tiny{Hopf}}}(H)$ of Hopf algebra automorphisms (this is just  $\Aut(H)$ if one thinks of $H$ 
as a bialgebroid over $\mathbb{C}$). 

We have seen that for a Galois object $A$ of a Hopf algebra $H$, the  
bialgebroid $\C(A, H)$ is a Hopf algebra. Also, the group of gauge transformations of the Galois object,
which is the same as the group of bisections $B(\C(A, H)$, is the group of characters of $\C(A, H)$ (see \eqref{isogen}). 
It turns out that these groups are also isomorphic to $\mbox{CoInn}(\C(A, H))$. 
If $\phi\in \mbox{Char}(\C(A, H))$, substituting $\phi^{-1}=\phi\circ S_{\C}$ in \eqref{altad}, 
for $X=a\ot a'\in \C(A, H)$, we get 
\begin{align}\label{autohopf1}
Ad_\phi (X) = \phi(\one{X}) \, \two{X} \, (\phi\circ S_{\C})(\three{X} ) = \mbox{coinn}(\phi)(X). 
\end{align}

As a particular case, let us consider again the Taft algebra $T_{N}$. We know from 
Sect.~\ref{taftal} that for any $T_{N}$-Galois object $A_s$ the bialgebroid $\C(T_{N}, A_{s})$ 
is isomorphic to $T_{N}$. And bisections of $\C(T_{N}, A_{s})$ are the same as characters of $T_{N}$ the group of which is isomorphic to $\IZ_N$. A generic character is a map $\phi_{r}:T_{N}\to \mathbb{C}$, 
given on generators $x$ and $g$ by 
$\phi_{r}(x)=0$ and $\phi_{r}(g)=r$ for $r$ a $N$-root of unity $r^N=1$.  
The corresponding automorphism $Ad_{\phi_{r}}=\mbox{coinn}(\phi_{r})$ is easily found to be on generators given by 
$$
\mbox{coinn}(\phi_{r})(g) = g , \qquad \mbox{coinn}(\phi_{r})(x) = r^{-1} x \, .
$$
It is known \cite[Lem. 2.1]{schau1} that
$\Aut(T_{N}) \simeq \Aut_{\mbox{\tiny{Hopf}}}(T_{N}) \simeq \mathbb{C}^{\times}$:  Indeed, given $r\in \mathbb{C}^{\times}$, one defines an automorphism $F_{r}: T_{N}\to T_{N}$ by $F_{r}(x):=rx$ and $F_{r}(g):=g$. Thus 
$Ad: \mbox{Char}(T_{N})\to \Aut(T_{N})$ is the injection sending $\phi_{r}$ to $F_{r^{-1}}$. 
Moreover, for $F\in \Aut(T_{N})$ and $\phi\in \mbox{Char}(T_{N})$, one checks that 
$Ad_{F\triangleright\phi}(x)=Ad_\phi (x)$ and $Ad_{F\triangleright \phi} (g )=Ad_\phi (g)$. Thus, 
as a crossed module, the action of $\Aut(T_{N})$ on $\mbox{Char}(T_{N})$ is the identity and the crossed module 
$(\mbox{Char}(\C(T_{N}, A_{s})),\Aut(\C(T_{N}, A_{s})), Ad, \id)$ is  
isomorphic to $(\mathbb{Z}_{N}, \mathbb{C}^{\times}, j, \id)$, with inclusion $j:\mathbb{Z}_{N}\to \mathbb{C}^{\times}$ given by $j(r):=e^{-\ii 2r\pi/N}$ and $\mathbb{C}^{\times}$ acting trivially on $\mathbb{Z}_{N}$.
\end{exa}

\subsection{More general bisections}
In parallel with the crossed module structure on bialgebroid automorphisms and bisections, there is a similar structure  on extended bisections as in Definition \ref{def:gebis} and `extended bialgebroid automorphisms'. 
These are pairs $(\Phi, \varphi)$ with $\varphi : B \to B$ an algebra automorphism 
and $\Phi: \cL \to \cL$ an unital invertible linear map, not required in general to be an algebra map, satisfying the equivariance properties in \eqref{amoeba(ii)}, while \eqref{amoeba(i)} is replaced by the bimodule property: $\Phi(b\triangleright a\triangleleft \tilde{b}) = \varphi(b)\triangleright \Phi(a) \triangleleft\varphi(\tilde{b})$, for $b, \tilde{b} \in B$.  They form a group $\Aut^{ext}(\cL)$ by map composition.

We sketch the construction that goes in the lines of Theorem \ref{thm. crossed module1}. Given any bisection 
$\sigma\in \B^{ext}(\C(A, H))$, definition \eqref{altad} 
gives a map $Ad_\sigma: \C(A, H)\to \C(A, H)$ that we repeat,
\begin{align}\label{Ad1}
Ad_\sigma(a\ot \tilde{a}) = \sigma(\one{(a\ot \tilde{a})}) \triangleright 
\two{(a\ot \tilde{a})} \triangleleft \sigma^{-1}(\three{(a\ot \tilde{a})})  .
\end{align}
This still covers the identity of $B$, that is 
$\varepsilon \circ Ad_\sigma \circ s = \id_B$, but is not an algebra map in general;
it is an extended automorphism of $\C(A, H)$. Indeed, it satisfies the properties
\eqref{amoeba(i)} and \eqref{amoeba(ii)}, and for two extended bisections $\sigma$ and $\tau$ 
one shows as before that $Ad_\sigma \circ Ad_\tau = Ad_{\tau\ast \sigma}$, and thus 
$Ad_\sigma$ is invertible with inverse $Ad_{\sigma^{-1}}$.
Moreover, if $\Phi$ is an extended automorphism of $\C (A, H)$ with inverse 
$\Phi^{-1}$ the formula \eqref{aaobs} is an action of $\Phi$ on $\B^{ext}(\C(A, H))$, a group 
automorphism of $\B^{ext}(\C(A, H))$. One really needs only to check that $\Phi \triangleright\sigma$ is well defined as an extended bisection since the rest goes as in the previous section. 
And with similar computations as those of Lemmas \ref{d2} and \ref{d3} one shows that
$Ad_{\Phi \triangleright\sigma} = \Phi^{-1}\circ Ad_\sigma \circ \Phi$, 
for any extended automorphism $\Phi$ and any $\sigma\in \B^{ext}(\C(A, H))$, and that 
$Ad_\tau \triangleright\sigma=\tau\ast\sigma\ast\tau^{-1}$, with any $\sigma$, $\tau\in \B^{ext}(\C(A, H))$.

We sum up the above with an analogous of Theorem \ref{thm. crossed module1}:
\begin{prop}\label{thm. crossed module2}
Given a faithfully flat Hopf--Galois extension $B=A^{coH}\subseteq A$, let $\C (A, H)$ be the corresponding Ehresmann--Schauenburg bialgebroid.
Then there is a group morphism 
$Ad: \B^{ext}(\C(A, H))\to \Aut^{ext}(\C(A, H))$ and an action 
$\triangleright$ of $\Aut^{ext}(\C(A, H))$ on $\B^{ext}(\C(A, H))$ such that the groups of extended  
automorphisms 
$\Aut^{ext}(\C(A, H))$ and of extended bisections $\B^{ext}(\C(A, H))$ form a crossed module.
\end{prop}

\begin{exa}
Consider a $H$-Galois object $A$ with $\C(A, H)$ the corresponding bialgebroid, a Hopf algebra itself. 
Given an extended bisection $\sigma \in \B^{ext}(\C(A, H)) \simeq \mbox{Char}^{ext}(H)$, 
in parallel with \eqref{autohopf} and \eqref{autohopf1} one can use \eqref{Ad1} to define 
an extended coinner automorphism of $\C(A, H)$,
$$
\mbox{coinn}(\sigma)(X):= 
Ad_\sigma(X) = \sigma(\one{X}) \, \two{X} \, \sigma^{-1}(\three{X}) ,
$$
for any $ X = a\ot \tilde{a}\in \C(A, H)$. This reduces to \eqref{autohopf1} when $\sigma = \phi$ is a character.
\end{exa}

\begin{exa}
In Example \ref{crossedmodule1} we gave an abelian crossed module for the Taft algebras. 
The use of extended characters and extended automorphisms yields a non-abelian crossed module. 
As we know the bialgebroid $\C(A_{s}, T_{N})$ of a Galois object $A_s$ for the Taft algebra $T_{N}$ is isomorphic to $T_{N}$. Thus $\Aut^{ext}(\C(A_{s}, T_{N})) \simeq 
\Aut^{ext}(T_{N})$ is the group of unital invertible maps
$\Phi: T_{N} \to T_{N}$ such that $\Phi(\one{h})\ot \Phi(\two{h})=\one{\Phi(h)}\ot\two{\Phi(h)}$ for $h\in T_{N}$. 

Let us illustrate this for the case $N=2$. The coproduct of $T_{2}$ on the generators $x,g$ will then require the following condition for an automorphism $\Phi$:
\begin{align*}
\one{\Phi(g)}\ot\two{\Phi(g)} &= \Phi(g) \ot \Phi(g) \nn \\
\one{\Phi(x)}\ot\two{\Phi(x)} &= 1 \ot \Phi(x) + \Phi(x) \ot g \nn \\
\one{\Phi(xg)}\ot\two{\Phi(xg)} &= g \ot \Phi(xg) + \Phi(xg) \ot 1 \, . 
\end{align*}
A little algebra then shows that
\begin{align*}
\Phi(g) = g \, , \qquad 
\Phi(x) =  c\, ( g - 1 ) + a_2 \, x \, , \qquad
\Phi(xg) = b\, ( 1 - g ) + a_1 \, xg 
\end{align*}
for arbitrary parameters $b, c \in \IC$ and $a_1, a_2 \in \IC^{\times}$ (for $\Phi$ to be invertible). 
As in \eqref{autverT} we can represent $\Phi$ as a matrix:
\begin{align}\label{authopfT}
\Phi : \begin{pmatrix}
1\\ xg\\ g\\ x\\
\end{pmatrix} 
\quad \mapsto \quad
 \begin{pmatrix}
1 & 0 & 0 & 0\\
b & a_1 & -b & 0\\
0 & 0 & 1 & 0\\
-c & 0 & c & a_2\\
\end{pmatrix}
\begin{pmatrix}
1\\ xg\\ g\\ x\\
\end{pmatrix} .
\end{align}
One checks that matrices $M_\Phi$ of the form above form a group: 
$\Aut^{ext}(T_{N}) \simeq \Aut_{\mbox{\tiny{Hopf}}}(T_{N})$. 

Given $\sigma \in \mbox{Char}^{ext}(T_{2})$ we shall denote $\sigma_a = \sigma(a) \in \IC$ for $a\in \{1, x, g, xg\}$. For the convolution inverse  $\sigma^{-1}$, from the condition $\sigma\ast\sigma^{-1}=\varepsilon$ we get on the basis that 
$$
\left \{
  \begin{aligned}
    &\sigma_1 = (\sigma^{-1})_1 = 1 \\
    &\sigma_{g}(\sigma^{-1})_{g}=1  \\
    &\sigma_{g}(\sigma^{-1})_{x}+\sigma_{x}=0\\
    &\sigma_{g}(\sigma^{-1})_{xg}+\sigma_{xg}=0 
  \end{aligned} \right. \qquad \Rightarrow \qquad 
\left \{
  \begin{aligned}
    &(\sigma^{-1})_{g} = (\sigma_{g})^{-1}  \\
    &(\sigma^{-1})_{x} = -\sigma_{x} (\sigma_{g})^{-1} \\
    &(\sigma^{-1})_{xg} = -\sigma_{xg} (\sigma_{g})^{-1}  
  \end{aligned} \right. \,\, .
$$
Then, computing $Ad_\sigma(h) = \sigma(\one{h}) \two{h} \sigma^{-1}(\three{h})$ leads to
\begin{align}
\label{ad-taft}
Ad_\sigma \begin{pmatrix}
1\\ xg\\ g\\ x\\
\end{pmatrix} 
=
\begin{pmatrix}
1 & 0 & 0 & 0\\
\sigma_{xg} & \sigma_{g} & -\sigma_{xg} & 0\\
0 & 0 & 1 & 0\\
- \sigma_{x} (\sigma_{g})^{-1} & 0 & \sigma_{x} (\sigma_{g})^{-1} & (\sigma_{g})^{-1}\\
\end{pmatrix}.
\end{align}
We see that the matrix \eqref{ad-taft} is of the form \eqref{authopfT} with the restriction that $a_2=a_1^{-1}$ so that $Ad_\phi$ has determinant $1$. Clearly, the image of $\mbox{Char}^{ext}(T_{2})$ form a subgroup of 
$\Aut^{ext}(\C(A_{s}, T_{2})) \simeq \Aut^{ext}(T_{N})$.
Moreover, $Ad: \mbox{Char}^{ext}(T_{2})\to \Aut^{ext}(T_{2})$ is an injective map. Finally, the action $Ad_{\Phi\triangleright\sigma}$ is represented by the matrix product:
\begin{align*}
& M_{Ad_{\Phi\triangleright\sigma}} = M_{\Phi^{-1}} \, M_{Ad_{\sigma}} \, M_\Phi \\
\qquad & =
 \begin{pmatrix}
1 & 0 & 0 & 0\\
a_1^{-1} [ \sigma_{xg} +  b (\sigma_{g} - 1 )] & \sigma_{g} & - a_1^{-1} [ \sigma_{xg} +  b (\sigma_{g} - 1 )] & 0 \\ ~~\\ 
0 & 0 & 1 & 0 \\ ~~\\
- a_2^{-1} [ \sigma_x \, (\sigma_g)^{-1} +c ( (\sigma_g)^{-1} - 1)] & 0 & a_2^{-1} [ \sigma_x \, (\sigma_g)^{-1} +c ( (\sigma_g)^{-1} - 1)] 
& (\sigma_g)^{-1} \\
\end{pmatrix}  \, .
\end{align*}
We conclude that as a crossed module the action on $\mbox{Char}^{ext}(T_{2})$ is not trivial.
\end{exa}

\appendix
\section{The classical gauge groupoid}\label{gaugegroupoid}

We collect here some basic facts of the gauge groupoid associated with a principal bundle and of the corresponding group of bisections \cite{KirillMackenzie}. Let $\pi: P\to M$ be a principal bundle over the manifold $M$ with structure Lie group $G$. Consider the diagonal action $(u, v)g:=(ug, vg)$ of $G$ on $P \times P$; denote by $[u, v]$
the orbit of $(u, v)$ and by $\Omega= P\times_{G} P$ the collection of orbits. Then $\Omega$ is a groupoid over $M$, 
--- the \emph{gauge} or \emph{Ehresmann groupoid} of the principal bundle,  with source and target projections given by
\begin{align*}
    s([u, v]):=\pi(v), \qquad t([u, v]):=\pi(u).
\end{align*}
The object inclusion $M\to P\times_{G} P$ is  
$m \mapsto \id_m :=[u, u]$, 
for $m\in M$ and $u$ any element in $\pi^{-1}(m)$. And the partial multiplication $[u, v']\cdot [v, w]$, defined when $\pi(v')=\pi(v)$ is 
\begin{align*}
    [u, v]\cdot [v', w]=[u, w g],
\end{align*}
for the unique $g\in G$ such that $v=v' g$.
Here one is really using the classical translation map $t: P\times_{M} P\to G$,
$(ug ,u) \mapsto g$.
One can always choose representatives such that $v=v'$ and the multiplication is then simply $[u, v]\cdot [v, w]=[u, w]$.
The inverse is
$$
[u, v]^{-1} = [v, u] .
$$

A \textit{bisection} of the groupoid $\Omega$ is a map $\sigma: M\to \Omega$, which is right-inverse to the source projection,
$s\circ \sigma=\id_{M}$, and is such that $t\circ \sigma : M \to M$ is a diffeomorphism. The collection of bisections,
denoted $\B(\Omega)$, form a group: given two bisections $\sigma_{1}$ and $\sigma_{2}$ their multiplication is defined by
\begin{align*}
    \sigma_{1}\ast\sigma_{2}(m):=\sigma_{1} \big((t\circ \sigma_{2})(m)\big)\sigma_{2}(m), \qquad \mbox{for} \quad m\in M .
\end{align*}
The identity is the object inclusion $m \mapsto \id_m$, simply denoted $\id$, with inverse given by
\begin{align*}
    \sigma^{-1}(m)= \big(\sigma\big((t \circ\sigma)^{-1}(m)\big)\big)^{-1};
\end{align*}
here $(t \circ\sigma)^{-1}$ as a diffeomorphism of $M$ while the second inversion is the one in $\Omega$.

The subset $\B_{P/G}(\Omega)$ of `vertical' bisections, that is those bisections that
are right-inverse to the target projection as well, $t\circ \sigma=\id_{M}$, form a subgroup of $\B(\Omega)$.

It is a classical result \cite{KirillMackenzie} that there is a group isomorphism between
$\B(\Omega)$ and the group of principal ($G$-equivariant) bundle automorphisms of the principal bundle
$$
\mathrm{Aut}_{G}(P) :=\{\varphi:P\to P \, ; \,\, \varphi (pg) = \varphi(p)g  \} \, ,
$$
 while $\B_{P/G}(\Omega)$ is isomorphic to
the subgroup of gauge transformations, that is principal bundle
automorphisms which are vertical (project to the identity on the base space):
$$
\mathrm{Aut}_{P/G}(P) :=\{\varphi:P\to P \, ; \,\, \varphi (pg) = \varphi(p)g  \, , \, \, \pi( \varphi(p)) = \pi(p) \} .
$$

\vspace{.5cm}
	
\noindent
{\bf Acknowledgments:}
We are most grateful to two anonymous referees for the great job with a previous version of the paper, in particular for bringing to our attention references that had escaped our attention and which lead to a 
completely rewritten paper. 
We thank Chiara Pagani for many useful discussions and for pointing us to the paper \cite{Be00}.
GL was partially supported by INFN, Iniziativa Specifica GAST and INDAM-GNSAGA.

\end{document}